\documentclass[12pt]{amsart}
\usepackage{amsmath,amssymb,mathrsfs,euscript,color,stmaryrd,mathabx}

\usepackage{epigraph}
\usepackage{hyperref}
\hypersetup{
pdfpagemode=UseOutlines,      % UseOutlines, UseThumbs, None, FullScreen : agencement au démarrage
pdfstartview=Fit,             % Fit, FitH, FitB, FitBH : vue de la page au départ (pleine largeur...)
pdffitwindow=true,            % bool: Maximiser
pdfpagelayout=TwoColumnsRight,% SinglePage, TwoColumnsRight/Left, OneColumn : affichage des pages
pdftoolbar=true,              % bool: Affichage de la barre d'outils
pdfmenubar=true,              % bool: Affichage de la barre de menus
bookmarksopen=false,          % bool: Dépliement des signets
bookmarksnumbered=true,       % bool: Numérotation des signets
colorlinks=true,              % bool: Liens colorés
pdfauthor={Ton nom},          % Auteur
pdftitle={Titre PDF},         % Titre
pdfcreator=PDFLaTeX,          % 
pdfproducer=PDFLaTeX,         %
linkcolor=blue,               % Couleur des liens
urlcolor=blue,                % url
anchorcolor=black,            % du texte
citecolor=blue,               % Couleur de citation 
frenchlinks=true,             % bool: Utiliser des petites majuscules pour les liens, plutôt que de la couleur
pdfborder={0 0 0}             % Ne pas encadrer les liens
}
\usepackage{xcolor}

\newtheorem{Theorem}[equation]{Theorem}
\newtheorem{Corollary}[equation]{Corollary}
\newtheorem{Proposition}[equation]{Proposition}
\newtheorem{Lemma}[equation]{Lemma}
\newtheorem{Hypothesis}[equation]{Hypotheses}
\theoremstyle{definition}
\newtheorem{Remark}[equation]{Remark}

\newtheorem{Definition}[equation]{Definition}

\makeatletter
\def\section{\def\@secnumfont{\mdseries}%
  \@startsection{section}{1}%
  \z@{.7\linespacing\@plus\linespacing}{.5\linespacing}%
  {\normalfont\scshape\centering}}
\def\subsection{\def\@secnumfont{\normalfont\bfseries}%
 \@startsection{subsection}{2}%
 \z@{.5\linespacing\@plus.7\linespacing}{-.5em}%
 {\normalfont\bfseries}}
\makeatother

\numberwithin{equation}{section}

\def\bdots{\mathinner{\mkern1mu\raise1pt\hbox{.}\mkern2mu\raise4pt\hbox{.}
           \mkern2mu\raise7pt\vbox{\kern7pt\hbox{.}}\mkern1mu}}

\def\Fx{F^\times}

\def\oF{{\mathfrak o}_F}
\def\oFx{{\mathfrak o}_F^\times}
\def\pF{{\mathfrak p}_F}
\def\kF{{k_F}}

\def\oE{{{\mathfrak o}_E}}
\def\oEx{{\mathfrak o}_E^\times}
\def\pE{{\mathfrak p}_E}

\def\fA{{\mathfrak A}}

\def\fM{{\mathfrak M}}

\def\fH{{\mathfrak H}}
\def\fJ{{\mathfrak J}}

\def\fR{{\mathfrak R}}

\def\End{{\operatorname{End}}}

\def\Hom{{\operatorname{Hom}\mspace{0.5mu}}}
\def\dim{\operatorname{dim}}

\def\det{\operatorname{det}}

\def\Ind{\operatorname{Ind}}
\def\cInd{\operatorname{c-Ind}}

\def\Lie{\operatorname{Lie}}

\def\val{\operatorname{val}}

\def\Jord{\operatorname{Jord}}

\def\IJord{\operatorname{IJord}}

\def\Sp{\operatorname{Sp}}
\def\GSp{\operatorname{GSp}}
\def\GL{\operatorname{GL}}

\def\SL{\operatorname{SL}}

\def\St{\operatorname{St}}

\def\tr{\operatorname{tr}}

\def\b{{\beta}}

\def\Supp{\operatorname{Supp}}
\def\diag{\operatorname{diag}}

\def\m{\mathbf m \medspace}

\def\adj{\,^a\negmedspace\,}
\def\Adj{\,^A\negmedspace}

\def\shaun#1{\textcolor{black}{#1}}  
\def\oldblue#1{\textcolor{black}{#1}}

\begin{document}

\title{Simple cuspidal representations of symplectic groups: Langlands parameter}

\author{Corinne Blondel}
\address{CNRS -- IMJ--PRG, Universit{\'e} Paris Cit\'e, Case 7012, 75205 Paris Cedex 13, France.}
\email{corinne.blondel@imj-prg.fr}

\author{Guy Henniart}
\address{Universit{\'e} de Paris-Sud, Laboratoire de Math{\'e}matiques d'Orsay, Orsay Cedex, F-91405.}
\email{guy.henniart@universite-paris-saclay.fr}

\author{Shaun Stevens}
\address{School of Mathematics, University of East Anglia, Norwich Research Park, Norwich NR4 7TJ, United Kingdom}
\email{Shaun.Stevens@uea.ac.uk}

%\thanks{} 

\begin{abstract} 
Let $F$ be a \shaun{non-archimedean} local field of odd residual characteristic. We compute the Jordan set of a simple cuspidal representation of a symplectic group over $F$, using explicit computations of generators of the Hecke algebras of covers reflecting the parabolic induction under study. 
When $F$ is a $p$-adic field we obtain the Langlands parameter of the representation. 
\end{abstract}

\date{\today}
%\subjclass[2010]{Primary 11F70, 22E50; Secondary .}
%\keywords{}
\maketitle

\begingroup % start a TeX group
\color{blue}% or whatever color you wish to use
\setcounter{tocdepth}{1}  %%% Shaun's edit
\tableofcontents
\endgroup

%%%%%%%%%%%%%%%%%%%%%%%%%%%%%%%%%%%%%%%%%%%%%%%%%%%%%%%
\section*{Introduction}
%%%%%%%%%%%%%%%%%%%%%%%%%%%%%%%%%%%%%%%%%%%%%%%%%%%%%%%

\shaun{Let~$F$ be a non-archimedean local field of residual characteristic~$p$, and let~$G$ be the group $\Sp(2N,F)$. The local Langlands conjecture for~$G$ attaches to a cuspidal representation\footnote{In this paper representations are smooth representations on complex vector spaces and by ``cuspidal representation'' we mean ``irreducible cuspidal representation.''}~$\pi$ of~$G$ a parameter of a Galois nature, or equivalently an irreducible representation~$\Pi$ of $\GL(2N+1,F)$. When $F$ has characteristic zero, the conjecture was established by Arthur~\cite{Arthur}, and M\oe glin~\cite{Moeglin} has shown that~$\Pi$ can be determined via the reducibility points of certain parabolically induced representations involving~$\pi$.}

\shaun{The method presented here to achieve this assumes that~$p$ is odd, and uses types and covers \`a la Bushnell--Kutzko~\cite{BK1} to obtain the reducibility points. It} was tested with success on $\SL(2,F)$ as early as 2009, in a joint project of the three authors initiated in January 2009 in Vienna.   The initial goal  was a full description of the L-packets of $\Sp(4,F)$ containing cuspidal representations  by means of types and covers. In those years, say 2009 to 2011, we did quite a lot of computations and 
completed a nice table presenting all cuspidal representations  of $\Sp(4,F)$, as classified in \cite{BS}, with the size of their packet and their neighbours in it. Some computations were done, but the expected tediousness of the others   made us choose a more conceptual way that we eventually wrote down in \cite{BHS}. We will explain this more precisely in a moment, let us just say that nonetheless, we accepted the idea that sometimes tedious computations can be useful to produce exact results, and we decided that the case of simple cuspidals of symplectic groups alone deserved such a treatment, along with the necessary work. This is the object of the present paper. 

\medskip 

\shaun{So let %$F$ be a \shaun{non-archimedean} local field of odd characteristic and let 
$\pi$ be a cuspidal representation of our symplectic group $G$.} % over $F$. 
We  need first to recall the main result in \cite{BHS}. The {\em Jordan set } $\Jord(\pi)$ of $\pi$ is the 
(finite) set of pairs $(\sigma, s)$ made of a self-contragredient cuspidal representation $\sigma$ of a  group  $\GL(k,F)$ for some $k$, and a real number $s\ge 1$,  such that, viewing $\GL(k,F)\times G$ as a maximal Levi subgroup of a suitable symplectic group $H$, 
 the normalised parabolically induced representation of $\sigma |\det|^s \otimes \pi$ to $H$ is reducible. When $F$ has characteristic zero, 
M\oe glin has shown that the Jordan set of $\pi$ determines the Langlands parameter of $\pi$. 

Theoretically  $\Jord(\pi)$ can be computed using types and covers,  thanks to the results of Bushnell and Kutzko that transform  parabolic induction in the groups  into  induction of modules over Hecke algebras, from the Hecke algebra of a type for the inertial class of  $\sigma |\det|^s \otimes \pi$ 
to the Hecke algebra of a cover of this type in $H$ \cite{BK1}. First of all one can associate to the representation $\pi$ a finite family $\mathscr F_\pi$ of {\it self-dual simple characters} and show that if $(\sigma, s)$ belongs to $\Jord(\pi)$, 
then $\sigma$ contains a simple character in the endoclass of the square of 
an element of $\mathscr F_\pi$. Then, having thus restricted the quest, we study the  cover of a type for the inertial class of $\sigma |\det|^s \otimes \pi$  for such a $\sigma$: the Hecke algebra of this cover has two generators which satisfy a quadratic relation computable  in a finite Hecke algebra deduced from the situation. 
Here results of Lusztig in finite reductive groups come into play, and eventually lead to a full description, not of $\Jord(\pi)$  itself, but of  the {\it inertial Jordan set} of $\pi$, which is the multiset  $\IJord(\pi)=\{([\sigma], s) \mid (\sigma, s) \in \Jord(\pi)\}$ (where 
$[\sigma]$ designates the inertial class of $\sigma$).  Indeed, the knowledge of   the finite reductive groups built from the underlying stratum of the cover,  of the level zero part of the cuspidal type and of Lusztig's results (see \cite{LS}),  produces   with a reasonable amount of computations (in particular of some characters \shaun{with} trivial square 
coming from  the compact subgroups involved in the construction) the quadratic relations satisfied by the generators, and eventually the inertial Jordan set.  

\medskip 
This is to be compared with the method presented here, that leads to the exact Jordan set if we are willing to pay the price of possibly very long computations, on a case-by-case basis.  This explains quite plainly why we changed   path towards \cite{BHS}. Yet obtaining exact results is definitely a respectable goal, more easily attainable whenever we deal with intertwining operators in one-dimensional spaces, i.e. when  $H^1=J^1$ \shaun{(in the standard notation of simple characters etc.)}, which occurs for simple cuspidals. Actually the computation below equally applies  whenever we deal with a stratum for which the extension field $F[\beta]/F$ is totally ramified of maximal degree $2N$, and it should  apply for other classical groups as well:   it has been used successfully in unitary groups in \cite{BT}.

\medskip 
There are alternatives to the computations that follow, some  are described in   \cite{BHS}. Indeed, once  we know the inertial Jordan set of $\pi$ in 
{\it loc.cit.},  we may sometimes fully determine some  parts of the Jordan set itself by  working on the Galois side, see    \cite[\S 7]{BHS} and section~\ref{Langlands} below. The computations presented here may nonetheless be necessary in severe cases where 
the ambiguity cannot be solved.  

\bigskip 

In the first section we present the method used to find the exact elements of the Jordan set. It is essentially an elaboration on the fundamental commutative diagram of \cite{BK1} -- the heart of the theory of covers -- in the case of a maximal Levi subgroup in a classical group. This diagram translates parabolic induction from $P$ to $G$ into induction of Hecke algebra modules, relying on a uniquely defined homomorphism of algebras $t_P$. Roughly speaking, when inducing from a maximal parabolic in a classical group, this morphism $t_P$ sends a generator of the Hecke algebra
on the fixed Levi component  $M$ of $P$ to the product of  two generators, say $T_0$ and $T_1$,  of the Hecke algebra over $G$. This equality  amounts to  normalising the corresponding intertwining operators  consistently.  This normalisation, in turn, allows for pinpointing the self-dual representation with ``highest reducibility value'' in the inertial class of the inducing representation (Theorem \ref{Thmethod}). 

In the second section we recall the definition of simple cuspidal representations in a symplectic group and we fix the notation for the particular simple cuspidal $\pi$ of $\Sp(2N,F)$ the Jordan set of which  we want  to compute. In particular we describe the underlying simple character $\psi_\beta$. We know from~\cite{BHS} (among other sources!) that this Jordan set  is 
$$
 \Jord(\pi)=  \{(\epsilon_1, 1),  (\sigma,1)   \}   
$$   
where $\epsilon_1$ is a character of $\Fx$   with trivial square and $\sigma$ is 
a cuspidal representation of $\GL(2N,F)$ attached to the simple character $\psi_{2\beta}$.  (Proposition \ref{BHSfirst}). 

In the third section we compute the character $\epsilon_1$ and in the fourth the simple cuspidal representation $\sigma$, using Theorem \ref{Thmethod} and precise computations of the coefficients of the quadratic relations satisfied by the generators of the Hecke algebra of the cover. At the end of section~\ref{dimension2Nsection} we explain how the case of a general simple cuspidal of $\Sp(2N,F)$ is easily deduced from the particular case that we have studied and we state the general result (Theorem~\ref{JordanSummary}). 

In the final section we go from Jordan set to Langlands parameter when $F$ has characteristic zero, or whenever the known results on the local  Langlands correspondence for $\Sp(2N,F)$ allow for such a move. We also discuss how our present result for simple cuspidal representations of $\Sp(2N,F)$ can also be obtained on the basis of the inertial Jordan set produced in \cite{BHS} together with a result of Lapid giving the $\varepsilon$-factor at $\frac 1 2$, whereas  there are cuspidal representations for which  this additional information is not sufficient. 

\bigskip 
{\bf Acknowledgements.}     
The authors  take the opportunity to signal the  work \cite{Sa2008} of Gordan Savin,
that determined the Jordan set of generic level zero cuspidal representations of classical groups, a
reference which inadvertently was absent from our paper \cite{BHS}.  

They  wish to thank    the organisers of  the workshop {\it Langlands Program:
Number Theory and Representation Theory }  in 
Oaxaca, Mexico,  for inviting each of  them to give a talk there in  December  2022. This gathering  gave them  the necessary impulse  to finish writing this important step of a long overdue project. 

\shaun{The third author was supported by EPSRC grants EP/H00534X/1 and EP/V061739/1.}

\section{Framework and method}\label{method}

\subsection{Covers and parabolic induction}

We go back to the founding paper by Colin Bushnell and Philip Kutzko \cite{BK1}. 
(We   use a mild variant of 
\cite{BK1} as explained in \cite{BB},  since we normalize parabolic induction and we use right-modules over the Hecke algebras, defined without contragredients.)

We fix $F$  a \shaun{non-archimedean} local field of \shaun{odd} residual characteristic $p$, we fix $G$ the group of $F$-points of a reductive algebraic group defined over $F$ and write $\mathfrak R(G)$ for the category of smooth complex representations of $G$. From now on all representations will be implicitly smooth and complex. 

 We fix  $M$ a  Levi subgroup of $G$, $P$ a parabolic subgroup of $G$ \shaun{with} Levi factor $M$, $U$ the unipotent radical of $P$, and $U^-$ the unipotent radical of the parabolic subgroup $P^-$ opposed to $P$ with respect to $M$. We fix a cuspidal inertial class $\mathfrak s_M=[M , \sigma]_M$ in $M$, which is the set of all twists 
$\sigma \chi$ of the  irreducible cuspidal representation $\sigma$ of $M$ by an unramified character $\chi$ of $M$, and denote by  $\mathfrak s_G=[M , \sigma]_G$ the corresponding inertial class in $G$, containing all pairs $G$-conjugate to some   $(M, \sigma \chi)$. We consider the functor 
$\Ind_P^G$  of normalized parabolic induction from the Bernstein block 
$\mathfrak R^{\mathfrak s_M}(M)$ in $\mathfrak R(M)$ to the Bernstein block $\mathfrak R^{\mathfrak s_G}(G)$ in $\mathfrak R(G)$. 

Assume that we have a type $(J_M, \lambda_M)$ for $\mathfrak R^{\mathfrak s_M}(M)$: so $J_M$ is a compact open subgroup of $M$, $\lambda_M$ is an  irreducible representation of $J_M$, hence finite-dimensional, acting on a space $V_{\lambda_M}$,  and all representations in 
$\mathfrak R^{\mathfrak s_M}(M)$ are generated by their $J_M$-isotypic component of type $\lambda_M$.  Then by \cite[Theorem 4.3]{BK1}, forming the Hecke algebra 
$$
  \begin{aligned}
\mathcal H(M , \lambda_M ) = \{f: M \rightarrow \End(V_{\lambda_M}) \mid f &\text{ compactly supported and } \\
&\forall g \in M, \ \forall j,k \in J_M, \ f(jg k) = \lambda_M(j) f(g) \lambda_M(k)
\},  
\end{aligned}
$$
 we have an equivalence of categories 
$$\mathfrak R^{\mathfrak s_M}(M)  \quad     \stackrel{E_{ \lambda_M}}{\longrightarrow}   \quad  \text{Mod-}\mathcal H (M, \lambda_M),  \qquad  
E_{ \lambda_M} (\omega)=  \Hom_{J_M} (\lambda_M, \omega), 
$$
where the structure of right-$\mathcal H (M, \lambda_M)$-module 
on $ \Hom_{J_M} (\lambda_M, \omega)$  is given by
\begin{equation}\label{defmodule}
\phi\cdot f (w) = \int_M \omega(g^{-1}) \phi ( f(g) w) dg \quad
(f \in \mathcal H (M, \lambda_M) ,\  \phi \in \Hom_{J_M} (\lambda_M, \omega),\  w \in V_{\lambda_M}).
\end{equation}

We further assume that we have a $G$-cover $(J_G, \lambda_G)$ of $(J_M, \lambda_M)$: a similar pair in $G$ with an Iwahori factorization  $
J_G=(J_G\cap U^-) (J_G\cap M) (J_G\cap U)
$,  
with $ \lambda_G$ trivial on  $J_G\cap U^- $ and $J_G\cap U$, with  
$J_G\cap M =J_M$ and $(\lambda_G)_{|J_M}=\lambda_M$, 
and with a strong additional condition that  
 provides an explicit injective homomorphism of algebras 
$$ t_P:  \mathcal H(M, \lambda_M) \longrightarrow \mathcal H(G, \lambda_G)  
$$
see \cite[Definition 8.1]{BK1}.  Then \cite{BK1} culminates with the 
assertion that $(J_G, \lambda_G)$ is a type for \shaun{$\mathfrak R^{\mathfrak s_G}(G)$} 
\cite[Theorem 8.3]{BK1}  and with   the 
 following commutative diagram that   transforms parabolic induction from   $\mathfrak R^{\mathfrak s_M}(M)$ to  $\mathfrak R^{\mathfrak s_G}(G)$     into  module    induction   over   Hecke algebras \cite[Corollary 8.4]{BK1}:
\begin{equation}\label{comm}
\begin{array}{ccc} \mathfrak R^{\mathfrak s_G}(G) &  \stackrel{E_{\lambda_G}}{\longrightarrow}  & \text{Mod-}\mathcal H(G, \lambda_G)
\\
  \Ind_P^{G} \uparrow   &   &   \uparrow   (t_P)_\ast   
\\ 
\mathfrak R^{\mathfrak s_M}(M)  &   \stackrel{E_{\lambda_M}}{\longrightarrow} & \text{Mod-}\mathcal H(M, \lambda_M)
\end{array}
\end{equation}
where, given a right $\mathcal H(M, \lambda_M)$-module $Y$, 
the $\mathcal H(G, \lambda_G)$-module $(t_P)_\ast   (Y)$ is 
the module $\Hom_{\mathcal H(M, \lambda_M)}(\mathcal H(G, \lambda_G), Y)$. 

\subsection{The equivalence of categories for cuspidal blocks}

We focus on the functor $E_{\lambda_M}$. By definition, irreducible objects in 
$\mathfrak R^{\mathfrak s_M}(M)$  form a single orbit under the group $X(M)$ of unramified characters of $M$, acting through $(\omega\chi) (g)= \chi(g)\omega(g)$ ($\chi \in X(M)$,
$\omega  \in \mathfrak s_M$, $g\in M$). The underlying space 
$E_{ \lambda_M} (\omega\chi)=  \Hom_{J_M} (\lambda_M, \omega\chi)$ is the same as  $E_{ \lambda_M} (\omega)$ because $\omega$ and $\omega\chi$ coincide on $J_M$, but those two spaces  differ as modules over $ \mathcal H(M, \lambda_M)$. The group $X(M)$ also acts on  $ \mathcal H(M, \lambda_M)$
by $(\chi f )(g)= \chi(g) f(g)$,   the action of $f  \in \mathcal H(M, \lambda_M)$ on $E_{ \lambda_M} (\omega)$ is the action of $\chi f$ on 
$E_{ \lambda_M} (\omega\chi)$.

\oldblue{When $M$ is a maximal Levi subgroup of a classical group $G$ and $p$ is odd,} 
 cuspidal representations of $M$ are known to satisfy the  following  conditions,  
 slightly stronger than \cite[(5.5)]{BK1}:

\begin{Hypothesis}\label{varpi}
 The type $(J_M, \lambda_M)$ satisfies the following. 
\begin{enumerate}
\item The intertwining of $\lambda_M$ is contained in a  compact mod center subgroup $\hat J_M$ of $M$, 
 containing $J_M$ as its unique maximal compact subgroup. 
\item 
$\lambda_M$ extends to $\hat J_M$ and for any such extension  $\widehat{\lambda}_M$ the representation 
$\cInd_{\hat J_M}^M \widehat{\lambda}_M$ is irreducible and cuspidal.
\item  There exists an element  $ {\varPi_{J_M}}$ of $ M$ such that 
$ \hat J_M = \varPi_{J_M}^\mathbb Z J_M$ and $\varPi_{J_M}^\mathbb Z \cap J_M = \{1\}$.  
\end{enumerate}
\end{Hypothesis}
From now on we assume that Hypotheses \ref{varpi} hold. Then  the Hecke algebra $\mathcal H(M, \lambda_M)$ is commutative \cite[Proposition 5.6]{BK1}, its irreducible modules are one-dimensional, they identify with  characters.  More precisely 
$\mathcal H(M, \lambda_M)$ is supported on $\varPi_{J_M}^\mathbb Z J_M$ and 
 isomorphic to $\mathbb C[\Psi, \Psi^{-1}]$ where  
$\Psi$ has support ${\varPi_{J_M}}  J_M=J_M  {\varPi_{J_M}}  $. This element $\Psi$ is  unique up to a non-zero scalar and characterized by the intertwining operator $\Psi({\varPi_{J_M}})
\in \End(V_{\lambda_M})$. Furthermore, if we pick an extension  $\widehat{\lambda}_M$ of $\lambda_M$ as in (ii), then the restriction of  $\widehat{\lambda}_M$ to  a compact subset of $\varPi_{J_M}^\mathbb Z J_M$ clearly belongs to $\mathcal H(M, \lambda_M)$, in other words $\Psi({\varPi_{J_M}})$ is a scalar multiple of $\widehat{\lambda}_M({\varPi_{J_M}})$. We would rather think about this the other way around: 
the Hecke algebra $\mathcal H(M, \lambda_M)$ does not depend on a particular choice 
of extension of $\lambda_M$, so we  fix a normalization of its generator $\Psi$ in an independent way, i.e. we consider the \shaun{non-zero} intertwining operator $\Psi({\varPi_{J_M}})$ chosen once and for all; in turn the extensions of $\lambda_M$ can be thought of relatively to $\Psi({\varPi_{J_M}})$. We introduce the following notation: 

\begin{Definition} We fix a normalization of $\Psi$ through the choice of an 
intertwining operator $\Psi({\varPi_{J_M}})$. 
Let     $\omega=\cInd_{\hat J_M}^M \widehat{\lambda}_M$ be an irreducible cuspidal representation of $M$ belonging to 
$\mathfrak s_M$,  where $ \widehat{\lambda}_M$ is an extension of $ \lambda_M$.
We let  $\zeta(\omega) $ be the  scalar such that 
$$\widehat{\lambda}_M({\varPi_{J_M}})=\zeta(\omega) \  \Psi({\varPi_{J_M}}).$$   
\end{Definition}

We  observe more closely the bottom line of  diagram \eqref{comm}. The functor  $E_{\lambda_M}$ attaches to each cuspidal representation 
$\omega$ in $\mathfrak s_M$, a character $\widecheck{\omega}$  of $\mathcal H(M, \lambda_M)$ that   represents the action of the algebra on 
$E_{\lambda_M}(\omega)$. This character is uniquely determined by its value on  $\Psi$ that we   compute using  formula \eqref{defmodule}, with 
$\phi \in \Hom_{J_M} (\lambda_M, \omega),\  v \in V_{\lambda_M}$:  
$$
\phi\cdot \Psi (v) = \int_{J_M} \omega(\varPi_{J_M}^{-1})\omega(g^{-1}) \phi (  \lambda_M(g) \Psi({\varPi_{J_M}}) v) dg  
=  \text{vol}(J_M) \   \omega(\varPi_{J_M}^{-1}) \phi (    \Psi({\varPi_{J_M}}) v)
$$
Since we have  $\omega=\cInd_{\hat J_M}^M \widehat{\lambda}_M$,   
the action of $\omega(\varPi_{J_M}^{-1})$ stabilizes the image of 
$\phi$ and acts on it as $\widehat{\lambda}_M(\varPi_{J_M}^{-1})$, so that   $\phi$ actually belongs 
to $\Hom_{\hat J_M} (\widehat{\lambda}_M, \omega)$. We get, after fixing the Haar measure on $M$ giving $J_M$ volume $1$:
$$
\phi\cdot \Psi (v) =      \omega(\varPi_{J_M}^{-1}) \phi (   \zeta(\omega)^{-1} \widehat{\lambda}_M({\varPi_{J_M}}) v) = 
\zeta(\omega)^{-1}  \phi(v). 
$$
 
\begin{Proposition}\label{hatomega} 
Let   $\omega$ be an irreducible cuspidal representation of $M$ belonging to 
$\mathfrak s_M$ and write   $\omega=\cInd_{\hat J_M}^M \widehat{\lambda}_M$ for  some extension   $ \widehat{\lambda}_M$   of $ \lambda_M$, 
characterized by 
$$  \widehat{\lambda}_M({\varPi_{J_M}})=\zeta(\omega) \  \Psi({\varPi_{J_M}}).$$
The action of $\mathcal H(M, \lambda_M)$  on 
$E_{\lambda_M}(\omega)$ is given by the character\footnote{\oldblue{This notation is convenient in our context but should not be confused with the usual notation for the contragredient representation. We will not use the latter.}}
$\widecheck{\omega}$ defined by  
$$\widecheck{\omega}  (\Psi)=  \zeta(\omega)^{-1} .$$
\end{Proposition}

Twisting  $\omega=\cInd_{\hat J_M}^M \widehat{\lambda}_M$ by $\chi  \in X(M)$ amounts to twisting $\widehat{\lambda}_M$ by $\chi_{|\hat J_M}$, hence we have 
$\zeta(\omega \chi)= \chi({\varPi_{J_M}}) \ \zeta(\omega )$ and: 
\begin{equation}\label{twist}
\widecheck{(\omega \chi)} (\Psi)=  \chi({\varPi_{J_M}})^{-1} \    \widecheck{\omega } (\Psi). 
\end{equation}

\subsection{The normalization of $\Psi$} From now on, we restrict   to the case studied in \cite{MS},  when $G$ is a classical group and  $M$ is a maximal Levi subgroup of $G$. The Levi subgroup $M$ identifies with a  direct product $M=\GL(k,F) \times M_0$ where $M_0$ is a classical group of the same sort as $G$. Unramified characters of $M$ have the form 
$$
(g, m_0) \longrightarrow |\det g|^s  \qquad \text{ with } s \in \mathbb C, \  
g \in \GL(k, F), \   m_0 \in M_0 . 
$$
The representation $\sigma$ of $M$ decomposes as $\sigma = \tau \otimes \pi$ 
where $\tau$ is a cuspidal representation of $\GL(k, F)$ and  $\pi$  a cuspidal representation of $M_0$; the irreducible objects of $\mathfrak s_M$ are the 
$ \tau |\det|^s \otimes \pi$ with  $s \in \mathbb C$. The type in $M$ for $\mathfrak s_M$ has the form 
$(J_M= J \times J_0, \lambda_M=\lambda\otimes \lambda_0)$ where 
$(J_0, \lambda_0)$ is a type constructed by \shaun{the third author} for $\pi$ and 
 $(J, \lambda)$ is a Bushnell--Kutzko type for $\tau$. In particular, we have 
$\pi=\cInd_{J_0}^{M_0} \lambda_0$,  and 
Hypotheses  \ref{varpi} hold  for $(J, \lambda)$: 
there are a compact mod center subgroup 
$\widehat J$ of $\GL(k, F)$ and an extension $\widehat \lambda$  of $\lambda$ to 
$\widehat J$ such that $\tau= \cInd^{\GL(k,F)}_{\widehat J} \widehat \lambda$, and  
there is   an element $\varpi_E$ such that $\widehat J = \varpi_E^\mathbb Z \times J$ with $\varpi_E^\mathbb Z \cap J = \{1\}$
(see 
\cite[\S 6]{BK}: the construction of $(J, \lambda)$ involves a finite extension $E$ of $F$ 
inside $M_k(F)$ such that the intertwining of $\lambda$ is $E^\times J$; 
we choose  a uniformizing element $\varpi_E$ of $E$ and remark that  
the ramification index of $E$ is uniquely attached to $\pi$). 
Then Hypotheses \ref{varpi} hold for $(J_M, \lambda_M)$ with 
$\widehat J_M=\widehat J \times J_0$, with $\widehat \lambda_M=\widehat \lambda \otimes \lambda_0$ and with ${\varPi_{J_M}}=(\varpi_E, 1)$. We let 
$(J_G, \lambda_G)$ be the $G$-cover of $(J_M, \lambda_M)$ built in  \cite{MS}. 

\medskip 

We have $|N_G(M)/M|=2$ and we 
further assume (see also  \cite[\S 11]{BK1}) that  the elements of $N_G(M)$ normalize $\mathfrak s_M$, which means that for some $s \in \mathbb C$, the representation $\tau |\det|^s$ is   self-dual (that is, equivalent to its contragredient in the symplectic and orthogonal cases, or equivalent to the contragredient of the conjugate representation in the unitary case). Theorem 1.2 in \cite{MS} essentially gives the following. 

\begin{Proposition}\label{varpi2}
\begin{enumerate}  
\item There are two elements $s_0$ and $s_1$ of $N_G(M)\backslash M$, belonging to open compact subgroups of $G$,  
that normalize $(J_M, \lambda_M)$ and  satisfy $s_0s_1= {\varPi_{J_M}}$.  
\item The Hecke algebra $\mathcal H(G, \lambda_G)$ is 
a two-dimensional module over $\mathcal H(M, \lambda_M)$  generated, as an algebra, by elements $T_0$ and $T_1$ of respective supports $J_G s_0 J_G$ and $J_G s_1 J_G$.
\item  The generators $T_0$ and $T_1$ can be normalized  
to satisfy     quadratic relations of the following shape: 
$$(T_i+1) (T_i-q^{r_i})=0,  \quad i=0,1,  \text{ with }  r_0, r_1 \ge 0.$$  
\end{enumerate}
\end{Proposition}

Indeed  $T_0$ and $T_1$ are defined up to non-zero scalar by their support,  normalizing them means    normalizing their values at $s_0$ and $s_1$ respectively, that are intertwining operators in the space of $\lambda_M$. This in turn provides a normalization for $\Psi$: we  impose
\begin{equation}\label{normZ}
 T_0 T_1= t_P(\Psi).  
\end{equation}

We make a quick comment here: the generators come from the $2$-dimensional Hecke algebras of two finite reductive groups, hence the quadratic relations, that    depend on the level zero part of the cover 
$(J_G, \lambda_G)$. They can be obtained through computations 
{\it \`a la Lusztig} in those finite groups, as in \cite[\S 5]{BHS}, and computation of a {\it twisting character} \cite[3.12 Theorem]{BHS}. \\

Proposition \ref{varpi2} implies that $\mathcal H(G, \lambda_G)$ has four characters, counted with multiplicity: the value of a character at $T_i $  is $-1$ or $q^{r_i}$, $i=0,1$. 
Hence the characters of $\mathcal H(M, \lambda_M)$ that induce reducibly to $\mathcal H(G, \lambda_G)$ are exactly the restrictions through $t_P$ of those four characters, their values at $\Psi$ belong to 
\begin{equation}\label{fourvalues}
\{ 1, -q^{r_0}, -q^{r_1}, q^{r_0+r_1}\}. \end{equation}  

\medskip

Now we recall what we know about reducibility on the group side.  
 The inertial class $\mathfrak s_M$  contains exactly two self-dual representations, 
say $ \tau_a \otimes \pi$ and $ \tau_b \otimes \pi$. For each of those there is a unique non negative  {\bf real} number $s_a$ or $s_b$ such that, for $s \in \mathbb R$:

$$
\Ind_P^G \tau_a |\det|^s \otimes \pi \text { reduces } \iff s=\pm s_a
$$

and the same  for $\tau_b, s_b$. So the four irreducible representations in $\mathfrak s_M$ (counted with multiplicities) that do NOT induce irreducibly are 

$$\{ \tau_a |\det|^{-s_a} \otimes \pi ,  \  \tau_a |\det|^{s_a} \otimes \pi,
  \   \tau_b|\det|^{-s_b} \otimes \pi, \   \tau_b |\det|^{s_b} \otimes \pi \}.$$

By Proposition  \ref{hatomega} 
  they correspond to the following characters of  $\mathcal H(M, \lambda_M)$, 
given by their value at $\Psi$: 
$$
|\det \varpi_E|^{s_a} \widecheck{\tau_a  \otimes \pi}(\Psi), \  
|\det \varpi_E|^{-s_a} \widecheck{\tau_a  \otimes \pi}(\Psi),  \  
|\det \varpi_E|^{s_b} \widecheck{\tau_b  \otimes \pi}(\Psi),  \  
|\det \varpi_E|^{-s_b} \widecheck{\tau_b  \otimes \pi}(\Psi).   
$$ 
This set of four values is identical to \eqref{fourvalues}, we deduce that 
one of $\tau_a$, $\tau_b$, say $\tau_a$, satisfies 
$$|\det \varpi_E|^{s_a} \widecheck{\tau_a  \otimes \pi}(\Psi) = 1 
\text{ and }  |\det \varpi_E|^{-2  s_a}  = q^{r_0+r_1}
$$
and the other one satisfies 
$$|\det \varpi_E|^{s_b} \widecheck{\tau_b  \otimes \pi}(\Psi) = -q^{\inf(r_0,r_1)}
\text{ and }  |\det \varpi_E|^{-2  s_b}  = q^{|r_0-r_1|}. 
$$

 We  find 
the corresponding self-dual representations 
of $\GL(k,F)$ with Proposition  \ref{hatomega}: 
$$\begin{aligned}
\tau_a&= \cInd_{\widehat J}^{\GL(k,F)}  \widehat \lambda_a  \text{ with } 
\widehat \lambda_a(\varpi_E) \otimes I_{V_{\lambda_0}}= |\det \varpi_E|^{s_a} \  \Psi(  {\varPi_{J_M}}  ),\\ \tau_b&= \cInd_{\widehat J}^{\GL(k,F)}  \widehat \lambda_b  \text{ with } 
\widehat \lambda_b(\varpi_E) \otimes I_{V_{\lambda_0}}= -q^{-\inf(r_0,r_1)}|\det \varpi_E|^{s_b}  \   \Psi(  {\varPi_{J_M}}  ).
\end{aligned}$$

We write for convenience $\equiv$ for ``equal up to a positive scalar''. 
The last touch is done by coming back to $T_0, T_1$, since $\Psi$ has been normalized by \eqref{normZ}, which is equivalent, up to a positive scalar that we don't need, to $T_0(s_0) T_1(s_1)\equiv \Psi(  {\varPi_{J_M}}  )$ \cite{BB}. We get finally:  

\begin{Theorem}\label{Thmethod}  We fix a cuspidal representation of the classical group $M_0$ and a self-dual cuspidal inertial class $\mathfrak s_k$  in $\GL(k,F)$, hence a cuspidal inertial class in 
the Levi subgroup  $M=\GL(k,F) \times M_0$ of the classical group $G$.
We fix as above a type $(J_M, \lambda_M)$  for this inertial class and a cover $(J_G, \lambda_G)$.  
  We normalize the two generators  $T_0$ and $T_1$ of the Hecke algebra  $\mathcal H(G, \lambda_G)$  
as in Proposition \ref{varpi2}, namely so that they satisfy quadratic relations
 $$(T_i+1) (T_i-q^{r_i})=0,  \quad i=0,1,   \text{ with }  r_0, r_1 \ge 0. $$ 

 The self-dual cuspidal representation  in   $\mathfrak s_k$  with the highest reducibility value  is the  self-dual cuspidal representation 
$
 \cInd_{\widehat J}^{\GL(k,F)}  \widehat \lambda_a $ 
characterized by $$
\widehat \lambda_a(\varpi_E) \equiv  T_0(s_0) T_1(s_1)   .  
$$
The other self-dual cuspidal  representation   in   $\mathfrak s_k$ is  
$\cInd_{\widehat J}^{\GL(k,F)}  \widehat \lambda_b $ with 
$\widehat \lambda_b(\varpi_E) \equiv -  T_0(s_0) T_1(s_1)   $.  
\end{Theorem} 

What we mean by ``highest reducibility value'' is: the representation having reducibility at $s_a$, since $s_a\ge s_b$, with equality if and only if one of $r_0, r_1$ is $0$. 
We recall that our goal is to determine, for a given cuspidal $\pi$ of $M_0$, the finite set 
of pairs $(\tau, s)$ with $\tau$ a cuspidal representation of some $\GL(k,F)$ and 
$s \in \mathbb R$, $s  \ge 1$, such that the normalized induced representation of 
$\tau |\det|^s \otimes \pi$ reduces.  We explained in \cite{BHS} how to construct this set, except possibly for an ambiguity between $\tau_a$ and $\tau_b$, in our notations above, that in some cases we couldn't solve. Theorem \ref{Thmethod} 
gives a way to solve the ambiguity. 
If  one of $r_0, r_1$ is $0$ we have no ambiguity to solve: indeed $s_a=s_b$ so either
$( \tau_a \otimes \pi, s_a )$ and $( \tau_b \otimes \pi, s_b  )$ 
both belong to the   set or neither of them does. Otherwise $s_a > s_b$ and  we produce the unique representation with reducibility at $s_a$. 

\begin{Corollary}\label{anypositive}
Theorem \ref{Thmethod} holds if $T_0$ and $T_1$ are \shaun{only} normalized  in such a way that they satisfy quadratic relations $T_i^2= b_i T_i + c_i$ with 
$b_i \ge 0$ and $c_i >0$, for $i=0,1$.
\end{Corollary}
 
Indeed such normalizations differ from the previous one by positive constants, 
that will not change the result, \shaun{except when a coefficient $b_i$ is $0$,  in which case both self-dual representations in the inertial class have highest reducibility value.}

 In the next sections,
we give examples  in which the 
computation is made easier by the fact that the intertwining operators are scalars. 

%%%%%%%%%%%%%%%%%%%%%%%%%%%%%%%%%%%%%%%%%%%%%%%%%%%%%%%
\section{Simple cuspidals}\label{section2}

\subsection{Definitions and notation}\label{notation}
We start with the necessary notation to describe simple cuspidal representations of symplectic groups, defined by  Gross and Reeder \cite[\S 9.2]{GR}. 

We let $F$ be a \shaun{non-archimedean} local field of odd residual characteristic $p$, with ring of integers $\oF$, maximal ideal $\pF$, residual field $\kF$ of cardinality $q_F=q$. We write $x \mapsto \overline x$ for the natural quotient map $\oF \rightarrow \kF$, 
and 
$\val (x)$ for the valuation of an element $x$ in $F$, \shaun{normalized so that~$\val$ has image~$\mathbb Z$}.  We fix  an additive character $\psi: F \rightarrow \mathbb C^\times$ with conductor $\pF$. We also fix for convenience a uniformizing element $\varpi_F$ of $F$. 
We let $\tilde G = \GL(2N,F)$, with centre $\tilde Z \simeq F^\times$,  and $G=\Sp(2N,F)$, the subgroup of $\tilde G$ preserving the 
alternating form $h_{2N}$ on $F^{2N}$ given by:

$$
h_{2N}(\left(\begin{smallmatrix} x_1\cr . \cr .\cr x_N\cr x_{N+1}\cr .\cr . \cr x_{2N} 
\end{smallmatrix}\right), \left(\begin{smallmatrix} y_1\cr . \cr .\cr y_N\cr y_{N+1}\cr .\cr . \cr y_{2N} 
\end{smallmatrix}\right) )=  x_1 y_{2N}+ \cdots + x_N y_{N+1}  - x_{N+1} y_N  - \cdots 
- x_{2N}y_1. 
$$

The matrix of the form $h_{2N}$ written in $N\times N$ blocks is 
$\left(\begin{smallmatrix} 0 & J_N\cr -J_N & 0
\end{smallmatrix}\right)$ where 
$J_N$ is the \shaun{$N$-by-$N$} matrix with $1$'s on the antidiagonal and $0$'s elsewhere.   The adjoint of a $2N$ by $2N$ matrix written in $N\times N$ blocks as 
$\left(\begin{smallmatrix} A & B\cr C & D
\end{smallmatrix}\right)$ is  $\left(\begin{smallmatrix} D^T & - B^T\cr 
- C^T &  A^T
\end{smallmatrix}\right)$  where %$A^T= J_N   A^t J_N$, i.e.  
$A \mapsto  A^T$ is  the transposition with respect to the \emph{anti}diagonal.
%%double-checked triple or quadruple checked 

The standard Iwahori subgroup $\tilde I_{2N}$ of 
$\tilde G$ is the fixator of the strict lattice chain $\Sigma_{2N}$ in $F^{2N}$  consisting of the columns of the order 
$\mathfrak A_{2N} = \left(\begin{smallmatrix} \oF& \oF&\dots & \oF & \oF \cr
 \pF& \oF& \oF&\dots& \oF \cr
\vdots &\ddots&\ddots &\ddots  &\vdots \cr 
  \pF&\dots& \pF& \oF& \oF \cr
  \pF& \pF&\dots& \pF& \oF \end{smallmatrix}\right)$
and their $\varpi_F^\mathbb Z $-multiples. 

The Jacobson radical of $\mathfrak A_{2N} $ is 
 $\mathfrak P_{2N}= \left(\begin{smallmatrix} \pF& \oF&\dots & \oF & \oF \cr
 \pF& \pF& \oF&\dots& \oF \cr
\vdots &\ddots&\ddots &\ddots  &\vdots \cr 
  \pF&\dots& \pF& \pF& \oF \cr
  \pF& \pF&\dots& \pF& \pF \end{smallmatrix}\right)$ with 
$\mathfrak P_{2N}^2= \left(\begin{smallmatrix} \pF& \pF& \oF &\dots & \oF \cr
 \pF& \pF& \pF& \oF &\dots\cr
\vdots &\ddots&\ddots &\ddots  &\vdots \cr 
  \pF&\dots& \pF& \pF& \pF \cr
  \pF^2& \pF&\dots& \pF& \pF \end{smallmatrix}\right)$, giving rise to subgroups   
$\tilde I_{2N}(1) = 1+ \mathfrak P_{2N}$ and $\tilde I_{2N}(2) = 1+ \mathfrak P_{2N}^2$
  of $\tilde I_{2N}$.

The successive maps $\textrm{I}_{2N} + (x_{i,j} )   \mapsto  
(x_{i,j} )$  and $(x_{i,j} )   \mapsto  (\overline x_{1,2}, \overline x_{2,3}, \cdots , \overline x_{2N-1, 2N}, \overline{\varpi_F^{-1} x_{2N,1}} )$ induce   isomorphisms
$$\begin{aligned}
\tilde I_{2N}(1)/ \tilde I_{2N}(2) \overset{\simeq}\longrightarrow \mathfrak P_{2N}/ \mathfrak P_{2N}^2
&\overset{\simeq}\longrightarrow  k_F^{2N}  .  
\end{aligned}
$$

\medskip

Taking now the intersections with $G$ we get the standard Iwahori subgroup $I_{2N}$ of $G$, with two subgroups $I_{2N}(1)$ and $I_{2N}(2)$,  and an isomorphism: 
$$\begin{aligned}
  I_{2N}(1)/   I_{2N}(2) 
&\overset{\simeq}\longrightarrow  k_F^{N+1}  \\
 (x_{i,j} ) \quad  &\longmapsto   (\overline x_{1,2}, \cdots ,\overline x_{N-1,N},  \overline x_{N, N+1}, \overline{\varpi_F^{-1} x_{2N,1}} ).  
\end{aligned}
$$

The center of $G$ is $Z \simeq\{\pm 1\}$. 
The {\it affine generic characters} of \cite[\S 9.2]{GR} are those characters of $  Z   I_{2N}(1)$ whose restrictions to $   I_{2N}(1)$ have the form 
$$
(x_{i,j} )   \longmapsto  \psi\left(  \alpha_1 x_{1,2}+   \cdots + \alpha_{N-1} x_{N-1, N}+ \alpha_{N} x_{N,N+1}  + \alpha_{2N} x_{2N,1}   \right) 
$$
with $\val(\alpha_i)=0$ for $i=1, \cdots,  N,$  and $\val(\alpha_{2N})=-1$. They  compactly induce irreducibly to cuspidal  representations of $  G$ called 
{\it simple cuspidal} representations of $G$.

\subsection{Description in terms of strata}\label{strata}

The chain $\Sigma_{2N}$, of period $2N$, can be scaled and translated into a unique lattice sequence $\Lambda_{2N}$ in $F^{2N}$ of period $4N$ and duality invariant $d=1$,   the usual convention in \cite{S5}.  That is, for $k \in \mathbb Z$, the dual lattice 
$$
\Lambda_{2N}(k)^\sharp = \{ X \in F^{2N} \mid h(X, \Lambda_{2N}(k) ) \subseteq \pF \}
$$ 
is equal to $\Lambda_{2N}(1-k)$. Note that  $\Lambda_{2N}(0)=\Lambda_{2N}(1)= \left(\begin{smallmatrix} \oF \cr . \cr .\cr \oF\cr \pF\cr .\cr . \cr \pF  
\end{smallmatrix}\right)$ ($N$ entries $\oF$, $N$ entries $\pF$).

According to \cite[\S 2]{BK2}, 
the natural filtration of $\mathfrak A_{2N} = \mathfrak A_0(\Lambda_{2N})$ given for integers $r$  by 
$$
\mathfrak A_r(\Lambda_{2N})=  \{ \phi  \in \End(F^{2N} ) \mid \forall k \in \mathbb Z \  \  
\phi(\Lambda_{2N}(k))\subseteq \Lambda_{2N}(k+r) \}  
$$ 
satisfies  $\mathfrak A_r(\Lambda_{2N})  = \mathfrak A_{\lceil \frac r 2 \rceil}(\Sigma_{2N})$ 
and $\val_{\Lambda_{2N}}= 2 \val_{\Sigma_{2N}}$, so that actually  
$\tilde I_{2N}(1) = 1+ \mathfrak A_1(\Lambda_{2N})= 1+ \mathfrak A_2(\Lambda_{2N})$ and $\tilde I_{2N}(2) = 1+ \mathfrak A_3(\Lambda_{2N})$. 

We leave aside the classification of affine generic characters and work directly with one whose restriction to   $I_{2N}(1)$   has the form 
$$x \longmapsto  \psi_\beta(x) = \psi\circ \tr (\beta(x-1))$$ for an element $\beta$ in $ \Lie (\Sp(2N,F)$ such that $\val_{\Lambda_{2N}}(\beta)= -2$ and $\beta^{2N}=(-1)^N \varpi_F^{-1}$. 
In particular $E=F[\beta]$ is a totally ramified   extension of $F$  of maximal degree $2N$. 
Actually we fix $\beta$ in $\mathfrak A_{-2}(\Lambda_{2N})$ as follows: 
 $\beta = \left(\begin{smallmatrix} 0& 0&\dots &\dots&\dots& 0 &  \varpi_F^{-1} \cr
 -1& 0& 0&\dots&\dots&\dots& 0 \cr
0 &\ddots&\ddots &\ddots& & &\vdots \cr 
\vdots   &\ddots  &-1 & 0 &\ddots &&\vdots \cr 
 \vdots  &  & 0 & 1 &0 &&\vdots \cr 
 \vdots  & &  & \ddots& \ddots & \ddots &\vdots  \cr
  0&0&\dots&\dots &0 & 1& 0 \end{smallmatrix}\right)$  
(with $N$ entries $-1$ and $N-1$ entries $1$).  
The adjoint of $\beta $ is  $-\beta$ and for 
$x= (x_{i,j} )$ we have 
$$\tr (\beta (x-1)) = - x_{1,2}-   \cdots - x_{N-1, N}- x_{N,N+1} + x_{N+1,N+2} + \cdots + x_{2N-1,2N} +  \varpi_F^{-1} x_{2N,1}   .$$  
Viewing $\psi \circ \tr (\beta(x-1))$ as a character of $  I_{2N}(1)/   I_{2N}(2) $ we can equate  
$x_{1,2}$ to $-x_{2N-1,2N}$ and so on, getting
\begin{equation}\label{psibeta} 
\psi \circ \tr (\beta(x-1))  = \psi ( - 2x_{1,2}-   \cdots -2 x_{N-1, N}- x_{N,N+1}   +  \varpi_F^{-1} x_{2N,1} ). 
\end{equation}

We note that the lattice chain underlying $\Lambda_{2N}$ is the set of $\beta^i \Lambda_{2N}(k)$ for $i \in \mathbb Z$ and any fixed $k$,   and that the $\oE$-order 
$\mathfrak A_0(\Lambda_{2N}) \cap E$ is just  the maximal $\oE$-order $\oE$. 

Thus $(\Lambda_{2N},2, 0, \beta)$ is a simple and maximal stratum in $\Lie (\Sp(2N,F))$, to which we  apply the machinery in \cite{S5}. Actually  $\beta$ is minimal over $F$ and   we have by \cite[\S 3.1]{S1}  $$J^1(\beta, \Lambda_{2N})= H^1(\beta, \Lambda_{2N})=I_{2N}(1), \qquad  
J(\beta,\Lambda_{2N})= Z J^1(\beta, \Lambda_{2N}) = Z I_{2N}(1) .$$   Then  $\psi_\beta$ is the unique simple character in $\mathcal C(\beta, \Lambda_{2N})$. The underlying stratum $(\Lambda_{2N},2, 0, \beta)$ is simple and maximal (attached to the totally ramified  field extension  of maximal degree),  so we obtain   the following  from \cite[3.6, 4.4 Theorem]{BHS}.  

\oldblue{
\begin{Proposition}\label{BHSfirst}
For any character $\chi$ of  the center $Z \simeq\{\pm 1\}$ of $G$, we consider   the  beta-extension 
$\kappa= \chi \otimes \psi_\beta$ of $\psi_\beta$, a representation of  $Z I_{2N}(1)$,   and  the  simple cuspidal  representation 
  $\pi=\cInd_{Z I_{2N}(1)}^G \chi \otimes \psi_\beta$  of $G$. \\
The Jordan set of $\pi$ is 
$
\   \Jord(\pi)=  \{(\epsilon_1, 1),  (\sigma,1)   \}   \   
$   
where $\epsilon_1$ is a character of $\Fx$   with trivial square and $\sigma$ is 
a cuspidal representation of $\GL(2N,F)$ attached to the simple character $\psi_{2\beta}$. 
\end{Proposition} 
We will discuss in the last section (see Theorem~\ref{LanglandsTh})  
 the Langlands parameter of $ \pi$.
}

In the next section we compute the character 
$\epsilon_1$, viewed as a character of $F^\times= \GL(1,F)$. In section \ref{GLN} 
we compute the cuspidal representation $\sigma$  of $\GL(2N,F)$.  
\oldblue{
Both computations rely on Theorem \ref{Thmethod}. There is a slight difference between them: in section \ref{dimension1section}  we use first \cite[4.4 Theorem]{BHS} to determine the restriction of $\epsilon_1$ to $\oFx$,  based on a twisting character    \eqref{twistingchar} computed in \ref{2.1}, then we proceed to the computation of the coefficients $b_0$ and $b_1$ using this restriction; in section \ref{dimension2Nsection} we proceed directly to the computation of $b_0$ and $b_1$
 keeping the restriction to $\oFx$ of the central character of $\sigma$ as a parameter, the value of which then results from the computation. Both ways are 
possible, we chose to use both. 
}

\section{The quadratic or trivial character}\label{dimension1section}

\subsection{The inertial Jordan set relative to the trivial endoclass}\label{2.1}

The four quadratic or trivial characters of $\Fx=\GL(1,F)$ are self-dual cuspidal representations attached to the null stratum $((\pF^k)_{k \in \mathbb Z}, 1,1, 0)$, 
the trivial character of $H^1((\pF^k),0)= 1+\pF$ and a self-dual beta-extension of this character to 
$J((\pF^k), 0) = \oFx$, which is  a quadratic or trivial character $\tau$ of 
$\oFx$.   In \cite[\S 3.6]{BHS} we built a cover $(J_P, \lambda_P)$ in $\Sp(2N+2,F)$ of the type 
$(\oFx \times Z I_{2N}(1), \tau \otimes \kappa)$ in the Levi subgroup 
$\GL(1,F) \times \Sp(2N,F)$. We recall some features of this cover. 

In the notation of  \cite[\S 3]{BHS} we have $V= F^{2N}$ as described above, $X=F^{2N+2}$ with elements written in coordinates $\left( x_0,  x_1, \cdots   ,  x_{2N},   x_{2N+1}
 \right)^t$ and  alternating form $h_{2(N+1)}$, and with 
$W=F$  the subspace given by the first coordinate $x_0$ and 
$W^\ast$ given by the last coordinate $x_{2N+1}$. On the space 
$W \oplus W^\ast$ we take the unique lattice sequence $\Lambda_2$ built on 
$\left(\begin{smallmatrix} \oF \cr \oF
\end{smallmatrix}\right)$, $\left(\begin{smallmatrix} \oF \cr \pF
\end{smallmatrix}\right)$ and their scalar multiples, that has period $4N$ and duality invariant $1$, and on $X=(W \oplus W^\ast) \perp V  $ we take the direct sum $\Lambda=\Lambda_2 \oplus \Lambda_{2N}$.  The stratum underlying the cover is $(\Lambda, 2,0, 0 \oplus \beta)$ (meaning: $0$ in $\Lie(\SL(2,F)$ and 
$\beta$ in $\Lie(\Sp(V)$). 

We  form two lattice sequences $\mathfrak M_0$ and $\mathfrak M_1$ in $X$ as follows. The first one $\mathfrak M_0$ (resp. the second one $\mathfrak M_1$) is the direct sum of the   unique lattice sequence $\mathfrak m_0$ (resp. $\mathfrak m_1$)  built on 
$\left(\begin{smallmatrix} \oF \cr \oF
\end{smallmatrix}\right)$ (resp.  $\left(\begin{smallmatrix} \oF \cr \pF
\end{smallmatrix}\right)$)  and its scalar multiples, that has period $4N$ and duality invariant $1$, and   $\Lambda_{2N}$. 

For the record, we first describe the relevant finite groups in our situation. For $i=0,1$, the 
finite group $P(\Lambda_\oE)/ P^1(\Lambda_\oE)$, isomorphic to 
$\oFx \times \{\pm 1\}$, is a Levi factor of $P(\Lambda_\oE)/ P^1(\mathfrak M_{i,\oE})$,
which is a parabolic subgroup of $\mathcal G_i= P(\mathfrak M_{i,\oE})/ P^1(\mathfrak M_{i,\oE})\simeq \Sp(2,\kF) \times \{\pm 1\}$.  Hence the two-dimensional Hecke algebras that arise here are just algebras on $\SL(2,\kF)$. We recall that,  for a character 
$\sigma$ of $\kF^\times$ with trivial square, viewed as a character of the Levi subgroup 
$\kF^\times$ of  $\SL(2,\kF)$, the Hecke algebra 
$\mathscr H(\SL(2,\kF), \sigma)$ has a generator $T$ satisfying the following quadratic relation: 
\begin{equation}\label{sl2}\begin{aligned}
(T-1)(T+1)&= 0   \text{ if } \sigma \ne 1,  \\ 
(T-q)(T+1)&= 0   \text{ if } \sigma = 1. 
\end{aligned}\end{equation} 

Now we apply \cite[4.4 Theorem]{BHS}. We are actually dealing with that part of the inertial Jordan set 
of $\pi$ 
relative to the trivial endoclass. The theorem says that it is the $\delta$-twist 
of the inertial Jordan set of the trivial representation of the trivial group, for a well-identified character $\delta$ of $\kF^\times$ which we will address shortly. 

The  Jordan set of the trivial representation of the trivial group has \shaun{itself long been known}: it has one element, the pair $(\mathbf \iota, 1)$. 
Indeed  the unique self-dual character $\sigma$ of $\GL(1,F)$ such that the normalized induced representation $\Ind_B^{\SL(2,F)} \sigma |.|^s$ reduces for some $s\ge 1$ (where $B$ is the standard Borel subgroup of upper triangular matrices) is the trivial character $\mathbf \iota$, and then $s= 1$   \cite[Corollary 9.3.3]{Ca}.  
%% Casselman works with normalized induction

The character $\delta$ is given by \cite[4.3 Proposition]{BHS} and can be computed  through \cite[4.2 Lemma]{BHS}: as a character of 
$\oFx$, its value at $x \in \oFx$ is the signature of the natural left action of $x$ on
$ \mathfrak J^1_{\mathfrak M_1}  \cap \Hom_F(V , W )/ \mathfrak H^1_{\mathfrak M_1}  \cap \Hom_F(V, W)$. (With the convention of {\it loc.cit.}  the space  
$V^0$ is the trivial space, hence $V^{\vee 0}=V$.) Implicit here is the stratum 
$(\mathfrak M_1, 2, 0 , 0 \oplus \beta)$ where $-2$ is the valuation of 
$0\oplus \beta$ relative to the sequence $\mathfrak M_1$, equal to 
$\val_{\Lambda_{2N}}(\beta)$. 

We must come back to the definitions. We recall that the \emph{jumps} of a lattice sequence $\Sigma$ in a vector space $S$ are  those integers~$i$ such that~$\Sigma(i) \ne \Sigma(i+1)$. The set of jumps of~$\Sigma$ is also the image of~$S\setminus\{0\}$ by the valuation map attached to~$\Sigma$, given  for~$y \in S\setminus\{0\}$ by~$\val_\Sigma(y) = \max \{ k \in \mathbb Z \mid y \in \Sigma(k)\}$.

Our stratum in $X=(W \oplus W^\ast) \perp V   $ is 
$( \mathfrak M_1, 2, 0, 0 \oplus \beta  )$ so the easiest way is to follow 
\cite[\S 3.3]{S1}. We obtain the $\oF$-orders $\mathfrak H_{\mathfrak M_1}$ and $ \mathfrak J_{\mathfrak M_1}$ written in blocks in the decomposition  $(W \oplus W^\ast) \perp V$: 
\begin{equation}\label{HetJ} \begin{aligned}
\mathfrak H_{\mathfrak M_1} &= \left(\begin{matrix}  
\mathfrak H (0,\mathfrak m_1)   & \mathfrak a^{12}_{2}(\mathfrak M_1)\cr
 \mathfrak a^{21}_{2}(\mathfrak M_1)&  \mathfrak H (\beta,\Lambda_{2N})   
\end{matrix}\right), \qquad  
  \mathfrak J_{\mathfrak M_1} &= \left(\begin{matrix}  \mathfrak J (0,\mathfrak m_1)   & \mathfrak a^{12}_{1}(\mathfrak M_1)\cr
 \mathfrak a^{21}_{1}(\mathfrak M_1)&  \mathfrak J (\beta,\Lambda_{2N})   
\end{matrix}\right). 
\end{aligned}
\end{equation}
 We concentrate on the first line of the upper-right block that corresponds to $ \Hom_F(V, W)$. To compare it between $\mathfrak H$ and $\mathfrak J$ we have to describe the lattices explicitly. We check that 
$\mathfrak m_1(t)= \left(\begin{smallmatrix} \pF^{[\frac{t+2N-1}{4N}]} \cr \pF^{[\frac{t+6N-1}{4N}]}
\end{smallmatrix}\right)$ (period $4N$, constant on the interval
$[-2N+1, 2N]$), the set of jumps of  $\mathfrak m_1$ is $4N\mathbb Z$. For $t \in [-2N, 2N-1]$ the lattices $\Lambda_{2N}(t)$ are the columns of the order $\mathfrak A_{2N}$ from right to left, each repeated twice; the set of jumps of $\Lambda_{2N}$ is 
the set of odd integers. 

The condition for some $b \in \Hom_F(V, W)$ to belong to  $\mathfrak a^{12}_{1}(\mathfrak M_1)$ or $\mathfrak a^{12}_{2}(\mathfrak M_1)$  is the following:
$$\begin{aligned}
b \in  \mathfrak a^{12}_{1}(\mathfrak M_1)   &\iff  \forall t
\text{ odd, }     b\Lambda_{2N}(t) \subseteq
\mathfrak m_1(t+1) \cap W \iff   b\Lambda_{2N}(-2N+1) \subseteq \oF   
\\
b \in  \mathfrak a^{12}_{2}(\mathfrak M_1)  &\iff  \forall t \text{ odd, }     b\Lambda_{2N}(t) \subseteq
\mathfrak m_1(t+2) \cap  W \iff 
\left\{\begin{matrix}b\Lambda_{2N}(-2N+1) \subseteq \oF  \cr
b\Lambda_{2N}(2N-1) \subseteq \pF 
\end{matrix}\right. 
\end{aligned}$$
So the condition is: all entries of $b$ in $\oF$ for $\mathfrak a^{12}_{1}(\mathfrak M_1)$, the first entry in $\pF$ and the others in $\oF$ for $\mathfrak a^{12}_{2}(\mathfrak M_1)$. Using  \cite[3.11 Lemma]{BHS}  we conclude that 
{\em 
\begin{equation}\label{twistingchar}\begin{aligned}
 &\delta   \text{ is the quadratic character of }   \oFx, \text{ in other words:} \\
&\IJord(\pi, \mathbf 1)= ( [\epsilon_1], 1) \text{  where }  \epsilon_1 
\text{ is a quadratic ramified character of }  F^\times. 
\end{aligned}\end{equation}}

We remark that $\IJord(\pi, \mathbf 1)$  does not depend on the character $\chi$ of $Z$ such that $\pi=\cInd_{Z I_{2N}(1)}^G \chi \otimes \psi_\beta$.

\subsection{The Jordan set relative to the trivial endoclass}\label{2.2}
\oldblue{We apply the  results of the first section to  $M=\GL(1,F) \times \Sp(2N,F)$ and $P$  the parabolic subgroup of $G^+= \Sp(2N+2,F)$ stabilizing the flag 
$\{0\} \subset W \subset W\oplus V \subset X$. Let $\epsilon$ be  a quadratic ramified character of $  F^\times$. We are studying normalized parabolic induction from   $M$ to
$\Sp(2N+2,F)$, specifically we are investigating the reducibility of the following representation: 
 $$ I(\pi, \epsilon, s ) = \text{Ind}_P^{G^+} \  \epsilon |\ |^s \otimes \pi \quad  (s \in \mathbb C) .   $$
We have a type $(\oFx \times Z I_N(1), \delta \otimes \kappa)$ in $M$ for  
$\epsilon   \otimes \pi$ and a cover $(J_P, \lambda_P)$ of this type in $G^+$. We ease notation by calling the respective Hecke algebras of these types 
$\mathcal H_M=\mathcal H(M ,  \delta \otimes \kappa) $ and $\mathcal H_{G^+}=\mathcal H(G^+ , \lambda_P )$.}
 
\oldblue{We consider the generator   $\Psi$ of $\mathcal H_M= \mathbb C [\Psi, \Psi^{-1}]$ supported on the 
$\oFx \times Z I_N(1)$-double \shaun{coset} of  
$\varPi_{J_M} =  \left(\begin{smallmatrix} \varpi_F & 0 &0 \cr
0 & I_{2N} & 0 \cr 0&0& \varpi_F^{-1}  \end{smallmatrix}  \right) $. 
The value of $\Psi$ at $\varPi_{J_M}$ is a non-zero intertwining operator of 
$\delta \otimes \kappa$; it is unique up to scalar, we take it as the identity on the space of $\kappa$ and some non-zero scalar   on the space of $\delta$.}  

\medskip 
\oldblue{
We turn to $\mathcal H_{G^+}$. The normalizer of $M$ in $G^+$ is the union of two $M$-cosets, the trivial coset and the coset of $t_0=\left(\begin{smallmatrix} 0& 0 &1 \cr
0 & I_{2N} & 0 \cr -1&0& 0 \end{smallmatrix}  \right)$ and $t_1=\left(\begin{smallmatrix} 0& 0 &- \varpi_F^{-1} \cr
0 & I_{2N} & 0 \cr \varpi_F&0& 0 \end{smallmatrix}  \right)$. We check that 
$t_0$ belongs to $P(\mathfrak M_{0,\oE})$, that $t_1$  
 belongs to $P(\mathfrak M_{1,\oE})$ and that $t_0t_1= \varPi_{J_M}$.  The algebra $\mathcal H_{G^+}$ has two generators $\mathcal T_0$ and $\mathcal T_1$ of respective supports $J_P t_0 J_P$ and $J_P t_1 J_P$, images of the corresponding generators of the two Hecke algebras in $\SL(2, k_F)$ described in \S \ref{2.1}. 
In view of \eqref{sl2} and \cite[3.14 Proposition]{BHS},  the only possibility for a reducibility at some $s$ with real part $1$ (a fact in our case, once chosen a self-dual base point)  is that both generators satisfy the quadratic relation $(T-q)(T+1)=0$, which defines them uniquely. Hence  $\mathcal T_0(t_0)$ and $\mathcal T_1(t_1)$ are uniquely 
determined by the quadratic relations $(\mathcal T_0-q)(\mathcal T_0+1)=0$ and $(\mathcal T_1-q)(\mathcal T_1+1)=0$.
}

\oldblue{
By Theorem \ref{Thmethod} and Corollary \ref{anypositive}, the quadratic character $\epsilon_1$ in the Jordan set of $\pi$ is characterised by:  
\begin{equation}\label{charepsilon1}
\epsilon_1(\varpi_F)\equiv    \mathcal T_0(t_0)   \mathcal T_1(t_1)
\end{equation}
where $\equiv $ means ``equal up to positive constant''. 
}
 
\oldblue{
\subsection{Computation of the argument of $\mathcal T_0(t_0) \mathcal T_1(t_1)$} 
We proceed   to determine    the arguments of 
$\mathcal T_0(t_0) $ and $\mathcal T_1(t_1)$  providing  quadratic relations with positive coefficients. 
}  
We work this out following \cite[\S 1.d]{BB}, that applies mutatis mutandis provided $t_0$ and $t_1$ behave well with respect to the Iwahori decomposition of  $J_P$, which we check first.

We write $P=MU$ for the parabolic subgroup defined in the previous subsection, 
with $U$ the unipotent radical of $P$, and we write $P^-=M U^-$ for the opposite parabolic with respect to $M$. We write $J_\Lambda $ for $J(\Lambda, 0\oplus \beta)$ and so on. From \cite[\S 3.6]{BHS} we have 
$$
J_P= (H^1_\Lambda \cap U^-)  (\oFx \times Z I_N(1)) (J^1_\Lambda \cap U)
$$

\subsubsection{Some lattice computations}

As in \eqref{HetJ},  following \cite[\S 3.3]{S1}, for the  stratum $( \Lambda, 2, 0, 0 \oplus \beta  )$, we write in blocks in the decomposition  $(W \oplus W^\ast) \perp V$: 
\begin{equation}\label{ForLambda} \begin{aligned}  
  \mathfrak H_{\Lambda} &= \left(\begin{matrix}  \mathfrak H (0,\Lambda_2)   & \mathfrak a^{12}_{2}(\Lambda)\cr
 \mathfrak a^{21}_{2}(\Lambda)&  \mathfrak H (\beta,\Lambda_{2N})   
\end{matrix}\right), \quad    \mathfrak J_{\Lambda} = \left(\begin{matrix}  \mathfrak J (0,\Lambda_2)   & \mathfrak a^{12}_{1}(\Lambda)\cr
 \mathfrak a^{21}_{1}(\Lambda)&  \mathfrak J (\beta,\Lambda_{2N})   
\end{matrix}\right), \\  
t_0 &= \left(\begin{matrix}  \left[\begin{smallmatrix} 0&1  \cr -1&0 
\end{smallmatrix}\right]   & 0\cr
0& I_{2N}  
\end{matrix}\right) , \quad
t_1 = \left(\begin{matrix}  \left[\begin{smallmatrix} 0&- \varpi_F^{-1}  \cr \varpi_F&0 
\end{smallmatrix}\right]   & 0\cr
0& I_{2N}  
\end{matrix}\right) .
\end{aligned}
\end{equation}
We write further $\mathfrak a^{12}_{i}(\Lambda) = \left(\begin{matrix}
R^1(i) \cr R^2(i)
\end{matrix}\right)$ where $R^1(i)$, $ R^2(i)$ are lattices of row vectors in $F^{2N}$ and similarly $\mathfrak a^{21}_{i}(\Lambda) = \left(C^1(i) \ C^2(i)
 \right)$ with lattices of column vectors. Recalling that $ J^1 (\beta,\Lambda_{2N}) 
=H^1 (\beta,\Lambda_{2N}) =I_{2N}(1)$ and that $H^1 (0,\Lambda_2) =J^1 (0,\Lambda_2) = I_1$, we get:
$$\begin{aligned}
J_P\cap U &=  \left(\begin{matrix} 1 & R^1(1) & \oF \cr
0 & I_{2N} & C^2(1) \cr 0&0&1
\end{matrix}\right) , \qquad \qquad \qquad
J_P \cap U^- =  \left(\begin{matrix} 1 &  0 & 0\cr
C^1(2) & I_{2N} & 0 \cr \pF& R^2(2)&1
\end{matrix}\right) \\
t_0( J_P\cap U^- ) t_0^{-1}&=  \left(\begin{matrix} 1 & R^2(2) & \pF \cr
0 & I_{2N} & C^1(2) \cr 0&0&1
\end{matrix}\right) , \qquad \qquad
t_0 (J_P \cap U)  t_0^{-1} =  \left(\begin{matrix} 1 &  0 & 0\cr
C^2(1) & I_{2N} & 0 \cr \oF& R^1(1)&1
\end{matrix}\right) \\
t_1 (J_P\cap U^- ) t_1^{-1}&=   \left(\begin{matrix} 1 & \varpi_F^{-1}R^2(2) & \pF^{-1} \cr
0 & I_{2N} &\tiny{ \varpi_F^{-1} C^1(2)} \cr 0&0&1
\end{matrix}\right) ,  \ 
t_1 (J_P \cap U)  t_1^{-1} \! = \! \left(\begin{matrix} 1 &  0 & 0\cr
\varpi_F C^2(1) & I_{2N} & 0 \cr \pF^2& \varpi_F R^1(1)&1
\end{matrix}\right)
\end{aligned}$$

We have to describe the lattices explicitly. We have seen before that for $t \in [-2N, 2N-1]$ the lattices $\Lambda_{2N}(t)$ are the columns of the order $\mathfrak A_{2N}$ from right to left, each repeated twice; the set of jumps of $\Lambda_{2N}$ is 
the set of odd integers. Now
$\Lambda_2$ has period $4N$, has a constant value equal to 
$\left(\begin{smallmatrix} \oF  \cr \pF 
\end{smallmatrix}\right)$ on the interval
$[-N+1, N]$, and the set of jumps of  $\Lambda_2$ is $N +2N\mathbb Z$.  

Elements $B= \left(\begin{matrix}
B_1 \cr B_2
\end{matrix}\right)$ of  $\mathfrak a^{12}_{i}(\Lambda)$, $i=1,2$, must satisfy 
$B \Lambda_{2N}(t) \subset \Lambda_2(t+i) $ for all $t$, i.e. 
$$\begin{aligned}
\underline{i=1} \quad &B \Lambda_{2N}(-N)\subset \left(\begin{smallmatrix} \oF  \cr \pF 
\end{smallmatrix}\right)   \quad \text{ and }  \quad  B \Lambda_{2N}(N)\subset \left(\begin{smallmatrix} \pF  \cr \pF 
\end{smallmatrix}\right) ; 
\\
\underline{i=2} \quad &  B \Lambda_{2N}(-N-1)\subset \left(\begin{smallmatrix} \oF  \cr \pF 
\end{smallmatrix}\right)   \quad \text{ and }  \quad  B \Lambda_{2N}(N-1)\subset \left(\begin{smallmatrix} \pF  \cr \pF 
\end{smallmatrix}\right) . 
\end{aligned}
$$
The first remark concerns parity. Since the jumps of $\Lambda_{2N}$ occur at odd integers, we have $\Lambda_{2N}(N)=\Lambda_{2N}(N-1)$ if and only if $N$ is odd. Hence  
if $N$ is odd we have $\mathfrak a^{12}_{1}(\Lambda)= \mathfrak a^{12}_{2}(\Lambda)$. 

We look at the rows of $B$ focusing on $R^1(1)$ and $R^2(2)$ which appear in $J_P$ above:  
$$\begin{aligned}
B_1  \in R^1(1) &\iff   B_1 \Lambda_{2N}(-N) \subset \oF \text{ and } 
B_1 \Lambda_{2N}(N) \subset \pF  \\
B_2  \in R^2(2) &\iff   B_2 \Lambda_{2N}(-N-1) \subset \pF \text{ (and } 
B_2 \Lambda_{2N}(N-1) \subset \pF).
\end{aligned}
$$
In particular $\varpi_F R^1(1) \subset R^2(2) \subset R^1(1) \subset \varpi_F^{-1}R^2(2)$, 

 and by duality   $\varpi_F C^2(1) \subset  C^1(2) \subset C^2(1) \subset \varpi_F^{-1}C^1(2)$.

Finally:  
\begin{equation}\label{Iwahori}
\begin{aligned}
t_0( J_P\cap U^- ) t_0^{-1} \subset\ &J_P\cap U 
\subset t_1 (J_P\cap U^- ) t_1^{-1}, \\
t_1 (J_P \cap U)  t_1^{-1} \subset\ 
&J_P \cap U^- \subset 
t_0 (J_P \cap U)  t_0^{-1} .
\end{aligned}
\end{equation} 

From these inclusions, we draw  that $t_1$ satisfies exactly the conditions in \cite[(1.3)]{BB}, whereas for $t_0$ we will only need to exchange the \shaun{roles} of $U$ and $U^-$.
With this the computation in \cite[\S 1.d]{BB}
 applies: we get the coefficients of the quadratic relations $ T^2=b_0T +c_0 \mathcal I $ and $ T^2=b_1T +c_1 \mathcal I $,  satisfied respectively by 
$\mathcal T_0$ and $\mathcal T_1$, from  \cite[(1.4)]{BB}. In particular 
{\it loc.cit.} provides immediately: 
\begin{Lemma}\label{c0c1}
The   coefficients $c_0$ and $c_1$ are positive if and only if $ \mathcal T_0(t_0)\mathcal T_0(t_0^{-1})$ and 
 $\mathcal T_1(t_1)\mathcal T_1(t_1^{-1})$   are positive, or equivalently $\delta(-1) \mathcal T_0(t_0)^2$ and 
$\delta(-1) \mathcal T_1(t_1)^2$  are positive.
\end{Lemma}

 For the coefficients $b_0$ and $b_1$ the computation based on  \cite[(1.4)]{BB} is more involved.

\subsubsection{Computation of the coefficient  $b_1$} 

  We must compute
$$
b_1= \sum_{j \in  (J_P\cap U) \backslash \Gamma }   \mathcal T_1 (j)  \qquad 
 \text{ where } \Gamma =   t_1 (J_P\cap U^- ) t_1^{-1} \cap J_P t_1 J_P.
$$
We have $J_P t_1 J_P =  (J_P\cap U^- )  (\oFx \times Z I_{2N}(1)) t_1 (J_P\cap U^- )$. The decomposition of an element $x$ of $J_P t_1 J_P$ as a product 
$$ x= u^-  \left(\begin{smallmatrix} \lambda & 0 & 0\cr
0 & g &0 \cr 0&0&\lambda^{-1}
\end{smallmatrix}\right)  t_1 v^- \text{ with } u^-, v^- \in J_P\cap U^-, 
 \lambda \in \oFx,                 g \in Z I_{2N}(1),$$
 is unique and gives 
$$\mathcal T_1 (x) =  \delta(\lambda) (\chi\otimes \psi_\beta)(g)\mathcal T_1(t_1).$$
To compute $b_1$ we must 
  work out   the matrix product to obtain,  by identification,   a characterization of  $\Gamma$ 
as some set of matrices $ j=  \left(\begin{smallmatrix} 1 & B & z \cr
0 & I_{2N} &C \cr 0&0&1
\end{smallmatrix}\right)$, with 
$B \in  \varpi_F^{-1}R^2(2)$, 
$C \in  \varpi_F^{-1} C^1(2)$ and $z \in \pF^{-1}$, and additional conditions, 
and compute $\lambda$ and $g$ as functions of $B$, $C$, $z$. We want:
$$
\begin{aligned}
 \left(\begin{smallmatrix} 1 & 0 & 0 \cr
D_1 & I_{2N} &0 \cr Z_1&H_1&1
\end{smallmatrix}\right)
 \left(\begin{smallmatrix} \lambda & 0 & 0\cr
0 & g &0 \cr 0&0&\lambda^{-1}
\end{smallmatrix}\right)
 \left(\begin{smallmatrix} 0& 0 &- \varpi_F^{-1} \cr
0 & I_{2N} & 0 \cr \varpi_F&0& 0
\end{smallmatrix}\right)
 \left(\begin{smallmatrix} 1 & 0 & 0 \cr
D_2 & I_{2N} &0 \cr Z_2&H_2&1
\end{smallmatrix}\right) &= 
\left(\begin{smallmatrix} -\lambda \varpi_F^{-1} Z_2 & -\lambda \varpi_F^{-1} H_2 &  -\lambda \varpi_F^{-1}  \cr
 -\lambda \varpi_F^{-1}D_1 Z_2+ g D_2 & g -\lambda \varpi_F^{-1} D_1H_2 & -\lambda \varpi_F^{-1} D_1 \cr y &Y& -\lambda \varpi_F^{-1} Z_1
\end{smallmatrix}\right) \\
&=  \left(\begin{smallmatrix} 1 & B & z \cr
0 & I_{2N} &C \cr 0&0&1
\end{smallmatrix}\right). 
\end{aligned}
$$
The obvious condition is that $z$ must have  valuation $-1$, then we let 
$\lambda=  - z  \varpi_F $. Next: 
\begin{itemize}
\item $Z_1=Z_2= z^{-1} \in \varpi_F \oFx$; 
\item $H_2= z^{-1}B  \in  R^2(2)$ and $D_1= z^{-1}  C \in C^1(2)$; 
\item  $g=  I_{2N}- z^{-1}CB $. 
\end{itemize} 
We must check  $g$. Conditions on $B$ and $C$ are  
 $ \varpi_F C \oF \subset  \Lambda_{2N}(N+2) $ and  $B \Lambda_{2N}(-N-1) \subset \oF$, they are equivalent by duality. From the second condition, the  entries in $B$ are in $\oF$ except the last $k$ ones in $\pF^{-1}$, for some $k$ with 
$1 \le k < N$. We will show that $ \varpi_F CB $ belongs to $\mathfrak A_0(\Lambda_{2N})$  if and only if all the entries of $B$ belong to $\oF$ -- this will show that actually $ \varpi_F CB $ belongs to $\mathfrak A_1(\Lambda_{2N})$. 
 
Recall that  $\left(\begin{smallmatrix} 1 & B & z \cr
0 & I_{2N} &C \cr 0&0&1
\end{smallmatrix}\right)$ belongs to $\Sp(2N+2,F)$ if and only if, writing 
$x_1$ to  $x_{2N}$ for the entries of $B$, left to right, and  $c_1$ to  $c_{2N}$ 
 for the entries of $C$, top to bottom,   we have 
$c_i= x_{2N-i+1}$ for $1\le i \le N$ and  $c_i= -x_{2N-i+1}$ for $N+1\le i \le 2N$, which we will write as $C= B^{\tau}$, and $BC=0$.  Assume that one of the last $k$ entries of $B$, say $x_{2N-j+1}$,  has valuation $-1$, then    the $(j, 2N-j)$ entry of  $ \varpi_F CB $ 
has valuation $-1$, which proves our claim.  In particular, when $g$ belongs to 
$\mathfrak A_0(\Lambda_{2N})$,  it belongs to $I_{2N}(1)$ and 
$(\chi\otimes \psi_\beta) (g)= \psi \circ \tr (- \beta  z^{-1}CB)$. 

\medskip

We leave aside for the moment the checking of the other coefficients and get on to computing $b_1$, with the following facts:
$$ \begin{aligned}
&\Gamma = t_1 (J_P\cap U^- ) t_1^{-1} \cap J_P t_1 J_P 
= \left\{  \left(\begin{smallmatrix} 1 & B & \varpi_F^{-1} u \cr
0 & I_{2N} &C \cr 0&0&1
\end{smallmatrix}\right) \in \Sp(2N+2,F) \mid  u \in \oFx, B \in \oF^{2N} \right\},  \\  
 &\mathcal T_1 ( \left(\begin{smallmatrix} 1 & B & z \cr
0 & I_{2N} &C \cr 0&0&1
\end{smallmatrix}\right)) =  \delta(-u) \psi \circ \tr(- \beta u^{-1}\varpi_FCB )\mathcal T_1(t_1) \quad  \text{ for } \left(\begin{smallmatrix} 1 & B & \varpi_F^{-1} u \cr
0 & I_{2N} &C \cr 0&0&1
\end{smallmatrix}\right) \in \Gamma
					.\end{aligned}$$ 																

We continue with the explicit element $\beta$ given in \S\ref{strata}, so that 
$$\begin{aligned}
\psi \circ \tr (-\beta u^{-1}\varpi_FCB)  &= \psi( u^{-1}\varpi_F (  2c_1x_2+ \cdots +2 c_{N-1}x_ N+ c_Nx_{N+1}   -  \varpi_F^{-1} c_{2N}x_1 )) \\
&=\psi( u^{-1}   x_1^2 ). 
\end{aligned}
$$
We need $b_1$ up to a positive constant, which we write as $\equiv$: 
 $$\begin{aligned}
b_1&\equiv \mathcal T_1(t_1)  \sum_{u \in k_F^\times} \delta(-u) \sum_{x \in k_F}\psi( u^{-1}   x^2 ) \equiv \mathcal T_1(t_1) \delta(-1) G(\delta, \psi), 
\end{aligned}
$$
where $G(\delta, \psi)$ is the Gauss sum $\sum_{u \in k_F^\times} \delta(u)  \psi( u  )$,  known to be   the product of $q^{\frac 1 2}$ and a square root of  
$(-1)^{\frac{q-1}{2}}$, namely 
\begin{equation}\label{Gauss}
\xi(\delta, \psi)= \frac{G(\delta, \psi)}{|G(\delta, \psi)|}  ,   \qquad 
\xi(\delta, \psi)^2=(-1)^{\frac{q-1}{2}}. 
\end{equation}
\begin{Proposition}
The normalization of $ \mathcal T_1 $ such that the coefficients $b_1$ and $c_1$ of the quadratic relation that it satisfies  are positive is given, up to a positive scalar, by 
$$
\mathcal T_1(t_1) =  \xi(\delta, \psi). 
$$
\end{Proposition}
Indeed, with this normalization the   coefficient $c_1$ is also positive, as stated in Lemma \ref{c0c1}, which stipulated that, up to a positive constant, 
$\mathcal T_1(t_1) $ was a square root of $\delta(-1)$. The exact square root is 
specified by the Gauss sum $G(\delta, \psi)$. 

As for the last \shaun{checks}: 
\begin{itemize}
\item $ -\lambda \varpi_F^{-1}D_1 Z_2+ g D_2=0 \iff z^{-1}C +D_2- z^{-1}CB D_2=0 $, 

which holds since $D_2= - H_2^{\tau} = -z^{-1} B^\tau=- z^{-1} C$ and $BC=0$. 
\item We have $Y=H_1g +H_2= H_1- z^{-1}H_1CB+ z^{-1} B$.  
   Since $H_1^\tau = -D_1= - z^{-1}C=- z^{-1}B^\tau$ we have $H_1= - z^{-1}B$  and $Y=0$ follows.  Then $y= - z^{-1}+ H_1 g D_2 +z^{-1} $ is $0$ for the same reasons. 
\end{itemize}

\subsubsection{Computation of the coefficient $b_0$ } As announced it is done in the same way with the \shaun{roles of} $U$ and $U^-$ being exchanged. We just write down the relevant facts. 
$$
b_0= \sum_{j \in  (J_P\cap U^-) \backslash \Gamma' }   \mathcal T_1 (j)  \qquad 
 \text{ where } \Gamma' =   t_0 (J_P\cap U ) t_0^{-1} \cap J_P t_0 J_P
$$
We have $J_P t_0 J_P =  (J_P\cap U )  (\oFx \times Z I_{2N}(1)) t_0 (J_P\cap U )$, and
$$ \begin{aligned}
&\Gamma' 
= \left\{  \left(\begin{smallmatrix} 1 & 0 & 0 \cr
D & I_{2N} &0 \cr u&H&1
\end{smallmatrix}\right) \in \Sp(2N+2,F) \mid   u \in \oFx, H \in (\oF, \cdots, \oF, \pF, \cdots, \pF)= \pF^{N} \times \oF^{N} \right\},  \\  
 &\mathcal T_1 (\left(\begin{smallmatrix} 1 & 0 & 0 \cr
D & I_{2N} &0 \cr u&H&1
\end{smallmatrix}\right)) =  \delta(-u) \psi \circ \tr(- \beta u^{-1}DH )\mathcal T_0(t_0) \quad  \text{ for }\left(\begin{smallmatrix} 1 & 0 & 0 \cr
D & I_{2N} &0 \cr u&H&1
\end{smallmatrix}\right) \in \Gamma'
.\end{aligned}$$ 	
Now $D$ and $H$ are related by  $D = - H^\tau$ so that $
\psi \circ \tr (-\beta u^{-1}DH)  = \psi(- u^{-1}    d_N^2  ) 
$,  and 
$$\begin{aligned}
b_0&\equiv \mathcal T_0(t_0)  \sum_{u \in k_F^\times} \delta(-u) \sum_{x \in k_F}\psi( - u^{-1}   x^2 ) \equiv \mathcal T_0(t_0)   G(\delta, \psi)  . 
\end{aligned}
$$  
\begin{Proposition}
The normalization of $ \mathcal T_0 $ such that the coefficients $b_0$ and $c_0$ of the quadratic relation that it satisfies  are positive is given, up to a positive scalar, by 
$$
\mathcal T_0(t_0) =   \delta(-1) \xi(\delta, \psi). 
$$
\end{Proposition} 

\subsection{Conclusion}\label{resultGL1}
Putting together \eqref{charepsilon1}   and the last two Propositions 
we obtain  
\begin{equation}\label{epsilon1}
\epsilon_1(\varpi_F)=  (-1)^{\frac{q-1}{2}} \xi(\delta, \psi)^{-2}=1.   
\end{equation}
In other terms, the Jordan set of $ \pi=\cInd_{Z I_{2N}(1)}^G \chi \otimes \psi_{\beta}$ relative to the trivial endoclass is 
$(\epsilon_1, 1)$ where $\epsilon_1$ is the ramified quadratic character 
such that $\epsilon_1(\varpi_F)= 1$. In terms of $\beta$, from \S \ref{strata}  we   replace $\varpi_F^{-1}$ by 
$(-1)^N \beta^{2N}= (-1)^{N+1} N_{E/F}(\beta)$  and get: 
$\epsilon_1( (-1)^{N+1} N_{E/F}(\beta))= 1$, or 
$$\epsilon_1(   N_{E/F}(\beta))=   (-1)^{(N+1)\frac{q-1}{2}}.$$ 
We remark that   the result does not depend on $\chi$, and   conclude: 

\begin{Proposition}\label{Jordanpi}
 The Jordan set of $ \pi=\cInd_{Z I_{2N}(1)}^G \chi \otimes \psi_{\beta}$ relative to the trivial endoclass is 
$(\epsilon_1, 1)$ where 
\begin{itemize}
\item $\epsilon_1$ is the ramified quadratic character 
that is trivial on the norms of $F[\beta]$ if $\frac{q-1}{2}$ is even or if $N$ is odd;  
\item $\epsilon_1$ is the ramified quadratic character 
that is \shaun{non-trivial} on the norms of $F[\beta]$ if $N$ is even  and $\frac{q-1}{2}$ is odd.
\end{itemize}
\end{Proposition}

%%%%%%%%%%%%%%%%%%%%%%%%%%%%%%%%%%%%
%%%%%%%%%%%%%%%%%%%%%%%%%%%%%%%%%%%%

\section{The simple cuspidal of $\GL(2N,F)$}\label{GLN}\label{dimension2Nsection}

We try and apply the same method to determine the simple cuspidal of $\GL(2N,F)$ 
that gives a reducibility with real part  $1$. We know from \cite{BHS}
the simple character underlying this representation: the square of the self-dual simple character extending  $\psi_\beta$. For the level zero part, section 5 in \cite{BHS} would give  the result, but we don't use it here. We   compute 
the generators of the Hecke algebra in order to describe completely the simple cuspidal.  

\subsection{The simple character and the cover}
We start again with the symplectic space $(V,h)=( F^{2N}, h_{2N})$ from section 
\ref{section2}. We work in the symplectic space 
$X= V \oplus V \oplus V$ equipped with the following symplectic form: 
$$
\mathbf h (\left(\begin{smallmatrix} a \cr b \cr c\end{smallmatrix}\right),   
\left(\begin{smallmatrix} a' \cr b' \cr c'\end{smallmatrix}\right))= 
h(a,c')    + h(b,b') + h(c, a') \qquad  ( a,b,c,a',b',c' \in V).
$$
We let $W=\left\{ \left(\begin{smallmatrix} a \cr 0 \cr 0\end{smallmatrix}\right)\mid a\in V\right\} $   and $W^\ast=\left\{ \left(\begin{smallmatrix} 0 \cr 0 \cr c\end{smallmatrix}\right)\mid c\in V\right\} $, and we make the identification 
$V =\left\{ \left(\begin{smallmatrix} 0 \cr b \cr 0\end{smallmatrix}\right)\mid b\in V\right\} $: this is  the symplectic space on which our original group $G=\Sp(2N,F)$ operates. For an endomorphism $Z $ of $V$ we denote by $\adj  Z $ the adjoint endomorphism, as described in \S \ref{notation}. 
For an endomorphism $Z $ of $X$ we  denote by $Z \mapsto \Adj Z$  the  adjoint map with respect to $\mathbf h$. We have
$$
\aligned
&\Adj\left(\begin{smallmatrix} g_1 & & \cr & g_2 & \cr & & g_3 
\end{smallmatrix} \right) 
= \left(\begin{smallmatrix} \adj g_3 & & \cr & \adj g_2 & \cr & & \adj g_1 
\end{smallmatrix} \right) ; 
  \quad   \quad 
\Adj\left(\begin{smallmatrix} & & g_1   \cr & g_2 & \cr   g_3 && 
\end{smallmatrix} \right) 
= \left(\begin{smallmatrix} &&\adj g_1 \cr & \adj g_2 & \cr  \adj g_3 && 
\end{smallmatrix} \right) ;
\\
&\Adj\left(\begin{smallmatrix} I & X & Z \cr & I & Y \cr & & I
\end{smallmatrix} \right) 
= \left(\begin{smallmatrix} I & \adj Y & \adj Z  \cr & I & \adj X \cr & & I 
\end{smallmatrix} \right)  ;
  \quad   \quad 
\Adj\left(\begin{smallmatrix} I & & \cr X & I & \cr Z & Y & I 
\end{smallmatrix} \right) 
= \left(\begin{smallmatrix} I & &  \cr   \adj Y & I &  \cr  \adj Z & \adj X & I
\end{smallmatrix} \right). 
\endaligned
$$

 We let $H=\Sp(X)\simeq \Sp(6N,F)$ and we consider  the embedding 
$$\begin{aligned}
\GL(W)\times G  \   &\longrightarrow  \    H \\  
(x, g ) \quad  &\longmapsto \quad 
\m (x,g) =  \left(\begin{matrix} x & & \cr &  g & \cr & & \adj x^{-1}
\end{matrix} \right) ,    \qquad x\in \GL(W), \  g \in G .  
\end{aligned}$$

 The image of $\m$ is 
a Levi subgroup 
$M$ of $H$. We let 
$P$ be  the parabolic subgroup of $H$ stabilizing the flag 
$\{0\} \subset W \subset W\oplus V \subset X$  and we write $P=MU$, 
with $U$ the unipotent radical of $P$, and    $P^-=M U^-$ for the opposite parabolic with respect to $M$.

 Each subspace $W$, $V$, $W^\ast$ of $X$ bears a  natural identification coordinate-wise with 
$F^{2N}$ through which we identify $\Lambda_{2N}$ to lattice sequences $\Lambda_W$, $\Lambda_V$, 
$\Lambda_{W^\ast}$. Note that $\Lambda_{W^\ast}$ is also 
the dual lattice sequence to $\Lambda_W$ when identifying 
 $W^\ast$ to the dual of $W$ through $\mathbf h $, i.e. 
 $$
 \Lambda_{W^\ast}(t)=\left\{ z\in W^\ast \mid  \forall x \in \Lambda_W(1-t)   \  \  
h(z,x) \in \pF \right\}.
$$

We recall our type  in $V$:
$$(J_V, \lambda_V)= (J(\beta, \Lambda_V), \chi\otimes 
\psi_{\beta}), \text{ with } J(\beta, \Lambda_V)= Z I_{2N}(1), $$   and consider the following data in $W$: 
\begin{itemize}
\item 
the simple and maximal stratum  
$(\Lambda_W,2, 0, 2\beta)$, 
\item the associated compact open subgroups 
$\tilde J^1(\beta, \Lambda_W)$ and $\tilde J(\beta, \Lambda_W)= \oFx \tilde J^1(\beta, \Lambda_W)$; 
\item  the simple character $\psi_{2\beta}$ of 
$\tilde J^1(\beta, \Lambda_W)$; 
\item a character $\delta$ of $\oFx $ with trivial square; 
\item  
the self-dual  type $ \  (\tilde J_W, \tilde \lambda_W)= (\tilde J(\beta, \Lambda_W), \delta\otimes 
\psi_{2\beta})  \   $ in $\GL(W)$. 
\end{itemize}
We form the 
 type $(J_M=\tilde J_W\times J_V, 
\lambda_M= \tilde \lambda_W\otimes \lambda_V)$ in $M$.

 We need a lattice sequence in $X$ which, together with 
$\beta_X=\beta \oplus \beta \oplus \beta$, will form a skew-simple stratum underlying 
an $H$-cover of  $(J_M, 
\lambda_M )$. The attached groups $H^1$, $J^1$ and $J$ must have  Iwahori decomposition with respect to $P=MU$. This will hold if 
the decomposition $X=W\oplus V \oplus W^\ast$ is properly subordinate to the stratum \cite[Corollaries 5.10, 5.11]{S5}, i.e.
\begin{itemize}
\item the lattices in the sequence are direct sums of lattices in $W,V,W^\ast$; 
\item  from one 
lattice to the next, at most one of the three parts changes. 
\end{itemize}
Using the definitions in 
 \cite[\S 2]{BK2}, we let 
$$
 \Lambda_X=(3\Lambda_W-2)\oplus (3\Lambda_V)\oplus (3\Lambda_{W^\ast}+2)$$
where $(3 \Lambda_W-2)(t)=3 \Lambda_W(t-2)=\Lambda_W([\frac{(t-2)+2}{3}])$, and so on.
The period of $\Lambda_X$ is $12N$. The dual of $\Lambda_X(t)$ is (with 
$1-[\frac x 3]=[\frac{1-x+4}{3}]$):
$$\begin{aligned}
\Lambda_W&(1-[\tfrac{(t+2)+2}{3}])\oplus \Lambda_V(1-[\tfrac{t+2}{3}])
\oplus \Lambda_{W^\ast}(1-[\tfrac{(t-2)+2}{3}])  \\& =
\Lambda_W([\tfrac{(1-t-2)+2}{3}])\oplus \Lambda_V([\tfrac{1-t+2}{3}])
\oplus \Lambda_{W^\ast}([\tfrac{(1-t+2)+2}{3}]) = \Lambda_X(1-t)
\end{aligned}$$
so $\Lambda_X$ has   duality invariant  $1$. The jumps of the sequence in $W$, resp. $V$, resp. $W^\ast$, occur for $t\equiv 5 $, resp. $t\equiv 3$, resp. $t\equiv 1$, mod $6$. We have $\Lambda_X(2t)=\Lambda_X(2t+1) $ for any $t \in \mathbb Z$, which implies  $\mathfrak A_{2t-1}(\Lambda_X)=\mathfrak A_{2t}(\Lambda_X)$ 
for $t \ge 1$.

We  form in $X$ the skew-simple stratum 
$
(\Lambda_X, 6, 0, \beta_X=\beta \oplus \beta \oplus \beta).
$ 
  We check the condition in \cite[\S 6.2]{S5}: the decomposition 
$X=W\oplus V \oplus W^\ast$ is exactly subordinate to the stratum, we have 
$\Lambda_X(1)=\Lambda_X(0)$ and  
$\Lambda_X(1)\cap  W^\ast \supsetneq \Lambda_X(2)\cap  W^\ast $. 
We  stick to the conventions and notations of {\it loc.cit.} and let 
$W=W^{(-1)}$, $W^\ast=W^{(1)}$, with $q_1=1$ and $q_{-1}=-1$; our parabolic subgroup $P$ is the same as in {\it loc.cit.}.

\bigskip  
We use  the cover    of 
$(J_M, \lambda_M)$ constructed by \shaun{the third author}  \cite[\S 6.2, \S 7.2.2]{S5}.  
 Since $\beta$ is minimal over $F$ and $\mathfrak A_{3}(\Lambda_X)=\mathfrak A_{4}(\Lambda_X)$, we have 
$J^1(\Lambda_X, \beta_X)= H^1(\Lambda_X, \beta_X)$ \cite[\S 3.1]{S1}.
  The skew-simple character $\psi_{\beta_X}$ of 
$J^1_X=H^1_X=H^1(\Lambda_X, \beta_X)$ restricts through $\m$  to  the character 
$\psi_{\beta_X}\circ \m = \psi_{2\beta}\otimes \psi_\beta$ of 
$ \tilde J^1(\beta, \Lambda_W)\times  J^1(\beta, \Lambda_{2N})$ and is trivial on the intersections with $U$ and $U^-$.  We have 
$$J_X:=J(\Lambda_X, \beta_X)=(H^1_X\cap U^-) \   \m\! ( \tilde J(\beta, \Lambda_W)\times  J(\beta, \Lambda_V)) \  ( H^1_X\cap U).$$ 
We get an $H$-cover $(J_X, \lambda_X)$ of $(J_M, \lambda_M)$ by letting $\lambda_X$ be trivial on $U$, $U^-$ and putting 
$\lambda_X\circ \m = \lambda_M$.

 \subsection{The Hecke algebra} 
We turn to $\mathcal H_{X}=\mathcal H(\Sp(X),\lambda_X) $. The normalizer of $M$ in $H$ is the union of two $M$-cosets, the trivial coset and the coset of the elements $s_1$ and $s_1^\varpi$   from  \cite[\S 6.2]{S5}: 
 $$
 s_1 = w_0= \left(\begin{matrix}
 0&0& I_{2N} \cr 0&I_{2N}&0 \cr 
 I_{2N} & 0 & 0 \end{matrix}\right) , \quad 
  s_1^\varpi =w_1=  \left(\begin{matrix}
 0&0&  \beta  \cr 0&I_{2N}&0 \cr 
 -\beta^{-1} & 0 & 0 \end{matrix}\right),
 $$
where we use    $\beta^{-1}$ as a uniformizing element for $E$, in other words 
we let $\beta^{-1}= \varpi_E$.   \\

In \cite[\S 7.2.2]{S5}, \shaun{the third author} constructs  self-dual lattice sequences  $\fM_0$ \shaun{and} $\fM_1$, of period $2$ over $E$, such that $w_0$ belongs to $P(\mathfrak M_{0,\oE})$ and  $w_1$  
 belongs to $P(\mathfrak M_{1,\oE})$. They are defined by 
$$
\fM_0 (2k+r)=  \left\{\begin{matrix} \varpi_E^k \Lambda_X(0) \text{ if } r=0,
\cr  \varpi_E^k \Lambda_X(1) \text{ if } r=1,
\end{matrix}\right.  
\qquad  
\fM_1 (2k+r)=  \left\{\begin{matrix} \varpi_E^k \Lambda_X(-2) \text{ if } r=0,
\cr  \varpi_E^k \Lambda_X(3) \text{ if } r=1.
\end{matrix}\right.
$$
 The algebra $\mathcal H_{X}$ has two generators $ T_0$ and $  T_1$ of respective supports $J_X w_0 J_X$ and $J_X w_1 J_X$.  Furthermore
$P_E(\Lambda_X)/P_E^1(\fM_i)$
is a maximal Levi subgroup of the finite  reductive group 
 $P_E(\fM_i)/P_E^1(\fM_i)$  and 
there is a quadratic character $\epsilon_{\fM_i}$ of $P_E(\Lambda_X)/P_E^1(\fM_i)$, depending only on 
$\fM_i$, $M$, $U$,   such that 
 $T_i$ satisfies a quadratic relation computed in 
$$\mathcal H(P(\mathfrak M_{i,\oE})/P^1(\mathfrak M_{i,\oE}),\epsilon_{\fM_i} (\delta\otimes \chi) ).$$

  Actually we are in the situation of  \cite[\S 3.16]{BHS}: 
the finite reductive groups obtained are $O(2,1)(\kF)$ and $\SL(2,\kF)\times \{\pm 1\}$. In the first one   the 
quadratic relation is always  \shaun{$T^2 = (q-1)T +q$}, the quotient of the roots is  $-q$ (i.e. $r_0=1$). In the second one, we get either the previous relation or \shaun{$T^2 = 1$}, the  quotient of the roots is $-q$ or $-1$ (i.e. $r_1=1 $ or $0$). Reducibility at $\pm 1$  corresponds to both relations equal to \shaun{$T^2 = (q-1)T +q$}.   We will come back to this later.

Now
$w_0$ and $w_1$ normalize  $J_X \cap M$ and exchange  $U$ and  $U^-$, and Lemma  7.11 in \cite{S5} 
gives 
$$
w_0 (J_X \cap U^- ) w_0^{-1} \subseteq J_X \cap U 
\qquad \text{ and } \qquad 
w_1 (J_X \cap U ) w_1^{-1} \subseteq J_X \cap U^-
$$
hence for $w_0$ : 
\begin{equation}\label{w_0}
\aligned 
&J_X w_o J_X = (J^1_X \cap U) w_0 J_M  (J^1_X \cap U), \\ 
&J_X \cap w_0 J_X w_0^{-1} = (J_X \cap U^- )J_M w_0 (J_X \cap U^- ) w_0^{-1}, 
\\
&\Omega_0:=J_X / J_X \cap w_0 J_X w_0^{-1} \simeq 
J_X \cap U / w_0 (J_X \cap U^- ) w_0^{-1} \simeq
J^1_X \cap U / w_0 (H^1_X \cap U^- ) w_0^{-1}; 
\endaligned 
\end{equation}
and for  $w_1$ : 
\begin{equation}\label{w}
\aligned 
&J_X w_1 J_X = (H^1_X \cap U^-) w_1 J_M   (H^1_X \cap U^-), \\ 
&J_X \cap w_1 J_X w_1^{-1} =  w_1 (J_X \cap U ) w_1^{-1}  J_M  (J_X \cap U ), 
\\
&\Omega_1:= J_X / J_X \cap w_1 J_X w_1^{-1} \simeq 
J_X \cap U^- / w_1 (J_X \cap U ) w_1^{-1} \simeq H^1_X \cap  U^- /  w_1 (J^1_X \cap U ) w_1^{-1} . 
\endaligned 
\end{equation}

We already know the possible forms of the quadratic relations satisfied by the generators, up to normalization. What we have to do is:  
\begin{enumerate}
	\item when two forms are possible, determine which one is obtained in terms of    $\chi$ and $\delta$; 
 \item determine, up to a positive scalar,  the normalization of the generators that 
gives a quadratic relation with positive coefficients -- in other words, choose the intertwining operator   $T_i(w_i)$ up to a positive scalar. 
 \end{enumerate}
Then   Theorem \ref{Thmethod} and the Corollary that follows will  give us the result.

We proceed, following the framework in  \cite[\S 1.d]{BB}. 
The  relations are     $T_i^2 = b_i T_i + c_i \mathbf 1$ where the scalars  $b_i$ and $c_i$ are given by the following formulae (simpler than in \cite{BB} since the space of   $\lambda_X$ has dimension $1$): 
\begin{equation}\label{bi}
\aligned 
&c_i  \ = \ |\Omega_i| \  T_i(w_i) \ T_i(w_i^{-1}) , 
\\
&b_i   \  =   \ \sum_{x \in \Omega_i} \ T_i(w_i^{-1} x^{-1} w_i)
=    \ \sum_{x \in Y_i} \ T_i(x),   
\endaligned
\end{equation}
where we let $Y_0 = (H^1_X \cap U^-) \backslash w_0^{-1} (J^1_X \cap U) w_0$ 
and   $Y_1 = (J^1_X \cap U) \backslash w_1^{-1} (H^1_X \cap U^-) w_1$. 

In the expression of   $b_i$, the support of the sum on $Y_i$ is the intersection 
of (a system of representatives of)   $Y_i$ with the  support of $T_i$. 
From the uniqueness of the Iwahori decomposition, the decomposition of some element as a product in  $U w_i M U$ or  
 $U^- w_i M U^-$ is unique (same reason: $P\cap U^- =\{1\}$). 
Let $x \in Y_0 \cap \Supp T_0$  and write  
$x=u w_o \shaun{d_0(x)} u^\prime$ with  $u, u^\prime \in J^1_X \cap U$ and 
 $\shaun{d_0(x)} \in J_M$, %we put  $\mathbf m = d_0(x)$ 
 and similarly for  $Y_1$ mutatis mutandis, 
consequently: 
\begin{equation}\label{bibi}
\aligned 
&b_i \ \  =  \ T_i(w_i)  
 \sum_{x \, \in \,  Y_i \cap \text{supp }T_i} \ \lambda_X(d_i(x)). 
\endaligned
\end{equation} 

\subsection{Relevant matrix  decompositions} 
We have to solve equations such as  
\begin{equation}\label{inf}
\left(\begin{matrix}
I &  &  \cr D & I &  \cr Z & H & I 
\end{matrix}\right) 
= \left(\begin{matrix}
I & B_1 & E_1 \cr  & I &  F_1\cr  &  & I 
\end{matrix}\right) 
\left(\begin{matrix}
 &  &  I \cr  & I &  \cr   I & & 
\end{matrix}\right) 
\left(\begin{matrix}
m &  &  \cr  & g &  \cr  &  &    \adj m^{-1}
\end{matrix}\right) 
\left(\begin{matrix}
I & B_2 &  E_2\cr  & I & F_2 \cr  &  & I 
\end{matrix}\right) 
\end{equation}
and 
\begin{equation}\label{sup}
\left(\begin{matrix}
I &  H &  Z \cr  & I & D  \cr & & I 
\end{matrix}\right) 
= \left(\begin{matrix}
I &  &  \cr F_1 & I &  \cr E_1 &  B_1 & I 
\end{matrix}\right) 
\left(\begin{matrix}
 &  &  \beta \cr  & I &  \cr  - \beta^{-1}  & & 
\end{matrix}\right) 
\left(\begin{matrix}
m &  &  \cr  & g &  \cr  &  &    \adj m^{-1}
\end{matrix}\right) 
\left(\begin{matrix}
I &  &  \cr F_2 & I &  \cr  E_2 & B_2 & I 
\end{matrix}\right) 
\end{equation}
in order to determine the intersections  $Y_i \cap \Supp T_i$, $i=0,1$.  
By uniqueness of the Iwahori decomposition, if the LHS belongs to the symplectic group, so do the elements in the RHS.  We want the LHS to belong to $Y_i$ and the elements in the RHS to belong to the relevant subgroups in the decomposition of   $J_X w_i J_X$, in particular we need $m \in \tilde J_W$, $g \in J_V$. 

(We remark that these equations are the ones considered by   Shahidi in \cite{Sh}, for orthogonal groups. They actually hold for $\GL(N')\times \Sp(2N)$ as well as the solutions below. Shahidi studies the relationship between  $m$ and $g$ in \eqref{inf},  
$m $ is almost $\adj Z^{-1}   $ or $Z$  and $g$ is related to the inverse of the ``norm'' of   $m$, namely   $-m^{-1} \adj  m$.)

We recall that  the adjoint of $\left(\begin{smallmatrix}
I &  &  \cr D & I &  \cr Z & H & I 
\end{smallmatrix}\right) $ is 
 $ \left(\begin{smallmatrix} I & &  \cr   \adj H & I &  \cr  \adj Z & \adj D & I
\end{smallmatrix} \right)$  so for such a matrix,   belonging to $\Sp(X)$ amounts to 
$H=- \adj D$ and $Z+ \adj Z+   \adj D D  =0$.  

To facilitate further \shaun{checks}, we \shaun{expand} the product on the RHS of \eqref{inf}: 
$$\left(\begin{matrix}
E_1 m  &T_1 m B_2+ B_1g & E_1 m E_2+ B_1gF_2+ \adj m^{-1}  \cr F_1 m  & F_1 m B_2 +g &  F_1m E_2+g F_2\cr m & m B_2   & m E_2 
\end{matrix}\right) $$

We see  that  \eqref{inf} has a solution if and only if  $Z$ is invertible,   given by 
\begin{equation}\label{solinf}
\aligned 
m &=   \, Z ; \quad 
B_2=- Z^{-1} \adj D ; \quad E_2 = Z^{-1} ; \quad F_1 = D Z^{-1} ; \quad E_1 = Z^{-1}  ;  \\
 g  &= I - (D Z^{-1}) Z ( - Z^{-1} \adj D )        = I + D Z^{-1}  \adj D .
\endaligned
\end{equation} As in 
 \cite[Corollary 3.2]{Sh} we have 
$g D  = D + D Z^{-1} \adj D D  = D - D Z^{-1}  (Z + \adj Z)= - D Z^{-1} \adj  Z $ so, when $D$ is invertible:
\begin{equation}\label{infDinversible}
g = - D Z^{-1}  \adj Z D^{-1} . 
\end{equation}

Similarly, the adjoint of $\left(\begin{smallmatrix}
I &  H &  Z \cr  & I & D  \cr & & I 
\end{smallmatrix}\right) $  is  $\left(\begin{smallmatrix}
I &  \adj D &  \adj Z\cr  & I & \adj H  \cr & & I 
\end{smallmatrix}\right) $
so belonging to  $\Sp(X)$ amounts to 
$H=- \adj D$ and $Z+ \adj Z  + \adj  D D=0$.  
The product on the RHS of \eqref{sup} is : 
$$\left(\begin{matrix}
\beta  \adj m^{-1} E_2  &\beta  \adj m^{-1} B_2 & \beta  \adj m^{-1}
\cr  g F_2  + F_1  \beta  \adj m^{-1}  E_2  &  g  + F_1  \beta  \adj m^{-1}  B_2&  F_1 \beta  \adj m^{-1}  \cr   -\beta^{-1} m + B_1 g F_2+T_1 \beta  \adj m^{-1}  E_2  & 
 B_1 g +E_1 \beta  \adj m^{-1}  B_2   & E_1 \beta  \adj m^{-1}
\end{matrix}\right) $$
so the general solution for \eqref{sup} is given, for an invertible  $Z $, by 
\begin{equation}\label{solsup}
\aligned 
 \adj m^{-1} &=  \beta^{-1} \,  Z ; \quad 
  B_2 = -Z^{-1}\adj D ; \quad E_2 = Z^{-1} ; \quad 
F_1= D Z^{-1}  ; \quad E_1 = Z^{-1} ;  
\\  g &= I - ( D Z^{-1})  Z  (-Z^{-1}\adj D)   = I + D Z^{-1} \adj  D  .
\endaligned
\end{equation}
Again $g D =  D + D Z^{-1} \adj  D D = D + D Z^{-1}(-Z- \adj Z)= - D Z^{-1}\adj Z$ and, 
when $D$ is  invertible: 
\begin{equation}\label{supDinversible}
g = - D Z^{-1} \adj Z D^{-1} . 
\end{equation} 
 
To proceed, we must describe the blocks in $J^1_X \cap U$ and other relevant subgroups. This is done in  \cite[Proposition 1]{Bl2} (for a lattice chain, but the lattice sequence $\Lambda_X$ is obtained by homothety-translation from the one in   
\cite{Bl2} and has same  $\fA_1$,  $\tilde H^1$ and $\tilde J^1$).  Here we have $t=3$ 
and a specially simple situation since 
$\fH^1= \fH^1(\b, \Lambda_{2N})= \fJ^1(\b, \Lambda_{2N})=\mathfrak A_{1}(\Lambda_{2N}):=\mathfrak A_{1}$. So:   

\begin{equation}\label{blocs}
\tilde J^1_X = \tilde H^1_X  = I + \left(\begin{matrix}
\mathfrak A_{1} & \oE+ \mathfrak A_{1} & \varpi_E^{-1} \mathfrak A_{1} \cr
\mathfrak A_{1} & \mathfrak A_{1} &  \oE+ \mathfrak A_{1}\cr
\pE + \varpi_E  \mathfrak A_{1} &\mathfrak A_{1} & \mathfrak A_{1} 
\end{matrix}\right) = I + \left(\begin{matrix}
\mathfrak A_{1} & \oF+ \mathfrak A_{1} &   \mathfrak A_{0} \cr
\mathfrak A_{1} & \mathfrak A_{1} &  \oF+ \mathfrak A_{1}\cr
\pE + \varpi_E  \mathfrak A_{1} &\mathfrak A_{1} & \mathfrak A_{1} 
\end{matrix}\right)  . 
\end{equation}

\subsection{Computation of $T_0$}
We are looking for solutions \eqref{solinf} of  \eqref{inf} such that   
\begin{itemize}
	\item 
   $\shaun{x=}\left(\begin{smallmatrix}
I &  &  \cr D & I &  \cr Z & H & I 
\end{smallmatrix}\right) 
$ \shaun{is in }%belongs to 
$w_0^{-1} (J^1_X \cap U) w_0$, 
i.e. $Z \in  \mathfrak A_{0}$ 
and $D \in \oF+ \mathfrak A_{1}$, modulo $H^1_X \cap U^-$;  
\item  \shaun{$x$} %the element 
belongs to $J_X w_0 J_X$, namely  
$B_1, B_2 \in  \oE+ \mathfrak A_{1}$, $E_1, E_2 \in  \mathfrak A_{0}$, 
$m \in \tilde J_W $ and $g \in J_V $. 
\end{itemize}
The first condition for existence is $Z \in \tilde J_W $. Then other constraints are obviously satisfied except the one for    $g$. But since  $\tilde J_W= \oFx +  \mathfrak A_{1}$, 
the condition  $Z+ \adj Z+ \adj D D=0$ implies $\adj D D \in  \oFx +  \mathfrak A_{1}$, 
which, added to $D \in \oF+ \mathfrak A_{1}$, implies $D \in \tilde J_W$. 
Then  
$g =- D Z^{-1}  \adj Z D^{-1} $ belongs to $\tilde J_W \cap \Sp(V)=J_V$.

We use \eqref{bibi} with notation in  \eqref{inf} and \eqref{solinf}. The general term in the sum is 
$$
\lambda_X(\shaun{d_0}(x))=   (\delta\otimes \psi_{2\beta}) (Z )\  
(\chi\otimes \psi_{\beta})(- D Z^{-1} \adj Z D^{-1}). 
$$
We write $Z=a(1+ z)$  and 
$D= u(1+d)$ with $a, u \in \oFx$ and $z, d \in  \mathfrak A_{1}$ and get: 
$$
\lambda_X(\shaun{d_0}(x))=  \delta(a) \psi_{2\beta} (1+z )\  
\chi(-1)  \psi_{\beta}( (1+d )(1+z)^{-1} (1+ \adj  z) (1+d)^{-1}) =   \delta(a)  
\chi(-1)   
$$ 
since $\psi \circ \tr (\beta \adj z)= \psi \circ \tr (- \beta  z)$.   
Now the sum in \eqref{bibi} is  on elements $Z,D \in \tilde J_W  $ with 
 $Z+ \adj Z+ \adj D D=0$, or equivalently on 
$a, u \in \oFx$, $d, z \, \in  \mathfrak A_{1}$, such that
$2a +  a(z + \adj  z ) + u^2  + u^2( d + \adj d ) + u^2 \adj d d = 0 $,  in particular $2a +  u^2\equiv 0 $ mod $\pF$.  
\shaun{Moreover, for each~$a,u$ satisfying this congruence, the number of pairs~$(d,z)$ satisfying the conditions is constant, independent of~$a,u$.} 
So, \shaun{working} up to a positive constant, we get

$$
\shaun{b_0    \equiv     T_0(w_0) }
 \sum_{\left.\begin{smallmatrix} a, u \in \kF^\times \cr 2a +  u^2 = 0  \end{smallmatrix}\right.
 } \ \delta(a)   
\chi(-1)  \equiv     T_0(w_0) \chi(-1) 
 \sum_{  u \in \kF^\times    } \ \delta(-u^2/2)     .  
$$ 
Since  $\delta$ is trivial on squares we have %up to positive constant:  
$  \shaun{b_0   \equiv
        T_0(w_0)  \chi(-1) \delta(-2)}
  $.    We know there is a   normalisation of 
 $T_0$ such that  $b_0= q-1$ and $c_0=q$. Since 
 $c_0  \ = \ |\Omega_0| \  T_0(w_0)^2$, this   normalisation  satisfies 
\begin{equation}\label{T0}
 \shaun{T_0(w_0)  \equiv   \chi(-1) \delta(-2) .} %\   C \quad \text{ with } C > 0 .  
\end{equation}

\subsection{Computation of $T_1$}
We look for solutions   \eqref{solsup} with: 
\begin{itemize}
	\item  $\shaun{x=}\left(\begin{smallmatrix}
I & H &  Z \cr   & I & D \cr   &   & I 
\end{smallmatrix}\right) 
$ \shaun{is in} %belongs to  
$  w_1^{-1} (H^1_X \cap U^-) w_1$, that is   $Z \in \beta (\oF+ \mathfrak A_{1})$ 
\shaun{and} $D \in \beta\fA_1$, mod $J^1_X \cap U$ ;  
\item  \shaun{$x$} %the element 
belongs to  $J_X w_1 J_X$, that is  
$B_1, B_2 \in \fA_1$, $E_1, E_2 \in \beta^{-1}  (\oF+ \mathfrak A_{1})$, 
$m \in \tilde J_W $ and $g \in J_V$. 
\end{itemize}
The first condition is   $m=- \beta \adj  Z^{-1} \in \tilde J_W$, that is   
$Z \in \beta\tilde J_W$. Then the other constraints are obviously satisfied, except the one for $g$ that we postpone.  We recall that  
$\varpi_E=\beta^{-1}$.  

The summation in \eqref{bibi} is  over the $(J^1_X \cap U)$-cosets of the intersection of $  w_1^{-1} (H^1_X \cap U^-) w_1$ with the support of $T_1$. 
An element of $w_1^{-1} (H^1_X \cap U^-) w_1$ can be written 
$\left(\begin{smallmatrix}
I & -\beta R  & - \beta U \beta \cr   & I & S \beta \cr & & I 
\end{smallmatrix}\right) $ with $U \in \pE + \varpi_E  \mathfrak A_{1}$ and 
$S \in \mathfrak A_{1} $, that is 
$\left(\begin{smallmatrix}
I &  H &  \beta z + t \cr  & I & D  \cr & & I 
\end{smallmatrix}\right) $  with  $z \in \oE$, $t \in \fA_0$, $D \in \fA_0$. The intersection with  $J_X w_1 J_X$ 
corresponds to  $z \in   \oEx$. We obtain a system of representatives of the quotient $Y_1$ as follows:  

 $\left(\begin{matrix}
I &  - \adj D&  \varpi_E^{-1}z - \frac 1 2 \adj D  D \cr  & I & D  \cr & & I 
\end{matrix}\right) $,  with  $z \in  k_F^\times$ and $D$  in  a system of representatives  $\fR$ that we detail later. We get 

$m =  \adj[\varpi_E (\varpi_E^{-1}z - \frac 1 2  \adj D D)]^{-1} =    \, z^{-1} \,  
\adj (1 - \frac 1 2 \varpi_Ez^{-1} \adj D  D)^{-1} $ and 

$g =  1 +  D(1 - \frac 1 2 \varpi_Ez^{-1} \adj  D D)^{-1} z^{-1}\varpi_E   \adj  D \equiv 1 +  D z^{-1}\varpi_E   \adj  D $ modulo 
$\fA_{3} $ (\shaun{recall that} $\Lambda_V$ has period $2$ over $E$).   
 
A term in the sum \eqref{bibi} can be computed as follows:   
$$
\aligned 
  &  \tilde \lambda_W( z^{-1} \,  
\adj (1 - \frac 1 2 \varpi_Ez^{-1} \adj D  D)^{-1} )\otimes \lambda_V(1 +  D  z^{-1}\varpi_E   \adj D ) \\
&=  \ \delta ( z^{-1}) \, \psi_{2\beta}(   \adj (1 - \frac 1 2 \varpi_Ez^{-1} \adj D  D)^{-1}  ) \,   \psi_\beta(1 +  D    z^{-1}\varpi_E  \adj  D)
\\
&=   \ \delta( z) \, \psi_{2\beta}(   
 1 - \frac 1 2 \varpi_Ez^{-1} \adj D D    ) \,  \psi_\beta(1 +  D  z^{-1}\varpi_E  \adj  D )
\\
&=   \delta( z) \,  \psi \circ \tr (2\b (- \frac 1 2 \varpi_Ez^{-1} \adj D  D)) \psi \circ \tr (\b D  z^{-1}\varpi_E   \adj D )
\\
&=   \ \delta( z) \,  \psi \circ \tr ( \b z^{-1} (-   \varpi_E \adj D  D  + D   \varpi_E    \adj D )). 
\endaligned 
$$ 

Remember that we took $\beta=\varpi_E^{-1}$ so 
$$
\aligned 
b_1 \,   \,  &=   \, T_1(w_1) 
 \sum_{D  \, \in \, \fR, \, z  \in k_F^\times
 } \ \delta( z) \,  \psi \circ \tr (  z^{-1} (-   \adj D  D  + \varpi_E^{-1} D   \varpi_E    \adj D )). 
\endaligned 
$$ 

Now   $ \fR$ is a system of representatives of  $ \fA_0 / \oE+\fA_1$, whereas for $D\in \oE$ the trace under $\psi$  is zero. We can use the bigger  quotient $ \fA_0 / \fA_1$ that has dimension $2N$ (see~\S\ref{notation}) and use for  $ \fR$ the diagonal matrices 
$D=\diag(d_1, \ldots, d_{2N})$ with coefficients in $\oF$ (mod $\pF$). 
Then 
$$
\aligned  &\adj D=\diag(d_{2N}, \ldots, d_1), 
\\& \varpi_E^{-1} D   \varpi_E   = \diag(d_{2N}, d_1,  \ldots, d_{2N-1}),    
\\& \tr ( \adj D  D  )= 2(d_1 d_{2N} + \cdots + d_N d_{N+1}), 
\\& \tr (  \varpi_E^{-1} D   \varpi_E   \adj  D )  = d_{2N}^2+ d_N^2+2(d_1 d_{2N-1} + \cdots + d_{N-1} d_{N+1}) . 
\endaligned 
$$ 
Working up to positive constant we get 
$$
\aligned 
b_1 \,   \,  &\equiv   \, T_1(w_1) 
 \sum_{d_1, \cdots, d_{2N}  \, \in \, \kF, \, z  \in k_F^\times
 } \ \delta( z) \,  \psi   (  z^{-1} ( d_{2N}^2+ d_N^2+2(d_1 d_{2N-1} + \cdots + d_{N-1} d_{N+1})  \\& \qquad \qquad \qquad \qquad \qquad \hskip5.5cm
-   2(d_1 d_{2N} + \cdots + d_N d_{N+1}))). 
\endaligned 
$$ 
Fixing all variables except one, say $d_k$ with  $k \ne N, 2N$,  we can factor out a sum 
$\sum_{d_k \in \kF} \psi  (u d_k)$, equal to $q$ if $u\in \pF$ and to 
$0$ if $\val(u)=0$.  So we are left with a sum with   conditions 
$d_{2N}=d_{2N-1}=\cdots = d_{N+1}$ and $d_N= d_{N-1}= \cdots = d_1$ and, always up to positive constant, to:  
$$
\aligned 
b_1 \,   \,  &\equiv  \, T_1(w_1) 
 \sum_{d_N, d_{2N}  \, \in \, \kF, \, z  \in k_F^\times
 } \ \delta( z) \,  \psi  (  z^{-1} ( d_{2N}^2+ d_N^2 
-   2 d_N d_{2N})  
\\ &\equiv     \, T_1(w_1) 
 \sum_{d_N, d_{2N}  \, \in \, \kF, \, z  \in k_F^\times
 } \ \delta( z) \,  \psi   (  z^{-1} ( d_{2N}- d_N)^2 
 )  \\ &\equiv     \, T_1(w_1) 
 \sum_{d   \, \in \, \kF, \, z  \in k_F^\times
 } \ \delta( z) \,  \psi   (  z^{-1} d^2 
 )  .  
\endaligned 
$$ 
 
If $\delta$ is trivial the sum over $z$ for a fixed $d$ is 
$q-1$ if $d=0$ and $-1$ if $d\ne  0$, so $b_1=0$.   Therefore 
we have reducibility at $1$ if and only if    $\delta $ is quadratic. If so, for a fixed $d$, 
the sum in $z$ is zero if    $d=0$, independent of  $d$ if $d$ is \shaun{non-zero}.   We obtain
 $$
\aligned 
b_1 \,  \,  &\equiv    \, T_1(w_1) 
 \sum_{  z  \in k_F^\times
 } \   \delta( z) \,  \psi   (  z )
 \equiv   \, T_1(w_1)   \xi(\delta, \psi ) .
\endaligned 
$$ 
where $\xi(\delta, \psi )$ is the normalised (modulus $1$)  Gauss sum defined in 
\eqref{Gauss}. 

We know that, if  $\delta $ is quadratic, there is a normalisation of 
 $T_1$ such that $b_1= q-1$ and  $c_1=q$. Since
 $c_1 \ = \ |\Omega_1| \  \delta(-   1  ) T_1(w_1)^2$ this  normalisation satisfies  
\begin{equation}\label{T1}
 \shaun{T_1(w_1) \equiv \xi(\delta, \psi )^{-1}.} % \   C \quad \text{ with } C > 0  .  
\end{equation}

\subsection{The answer}\label{answerGL2N}
 We fix now $\delta$ as  the \shaun{(non-trivial)} quadratic character of $\oFx$. 
The cuspidal type $(\tilde J_W, \tilde \lambda_W)$ extends to 
the compact mod center subgroup $E^\times \tilde J_W$ by choosing a character 
$\tau$ of 
$E^\times$ extending $\delta$. This is equivalent to choosing the value 
of $\tau $ on  a uniformizing element of $E$. The induced representation of $\tau \otimes \psi_{2\beta}$ to $\GL(W)$ is then cuspidal irreducible. 

There is exactly one of these representations, say $\sigma=\cInd_{E^\times \tilde J_W}^{\GL(2N,F)}  \tau\otimes \psi_{2\beta}  $,   such that $\sigma$ is self-dual and $\Ind_P^G \sigma |\det| \otimes \pi$  
is reducible. 
This representation is characterized by the value of $\tau $ on a uniformizing element  given by Theorem \ref{Thmethod}. Since we have 
$$
w_0 w_1 =  \left(\begin{matrix} -\beta^{-1} & 0 &0 \cr
0 & I_{2N} & 0 \cr 0&0& \beta  \end{matrix}  \right)
$$
we must have,   up to a positive constant: 
$$
\tau( -\beta^{-1})   \equiv    \chi(-1) \delta(-2)   \xi(\delta, \psi )^{-1} 
$$
But the representation must be self-dual and the inducing character also, hence
$$
\tau( \beta^{-1})   \equiv    \chi(-1) \delta(2)   \xi(\delta, \psi )^{-1} 
$$

\begin{Proposition}\label{JordanpiGL}
 The Jordan set of $ \pi=\cInd_{Z I_{2N}(1)}^G \chi \otimes \psi_{\beta}$ relative to the  endoclass of  the simple character $\psi_{2\beta}$ of $\tilde I_{2N}(1) $ is 
$\Jord(\pi, \psi_{2\beta})= \{ (\sigma, 1)\}$ with  
 $$
\sigma= \cInd_{E^\times \tilde I_{2N}(1)}^{\GL(2N,F)}  \tau\otimes \psi_{2\beta}  
$$
where $\tau_{|\oEx} $ is the quadratic character of $\oEx$ and 
$$\tau(\beta) =\chi(-1)\delta(2)  \xi(\delta, \psi ).$$
\end{Proposition}

We notice that $\tau(-\beta^2)= \tau(-1) \xi(\delta, \psi )^2=1$ and   that 
$$
\tau(-\beta^{2N})= \delta(-1) \xi(\delta, \psi )^{2N} = [(-1)^\frac{q-1}{2}]^{N+1}
$$
is  trivial if $N$ is odd and equal to $(-1)^\frac{q-1}{2}$ if $N$ is even.  

\subsection{Other simple cuspidals}\label{othersimple}

So far we have computed the Jordan sets of  the simple cuspidal representations of $G=\Sp(2N,F)$ whose restriction to $I_{2N}(1)$ is given by the element $\beta$  of~\S\ref{strata}. Note, however, that $\beta$  depends on the choice of the uniformizer  
$\varpi_F$ of $F$, which we had fixed but is otherwise arbitrary. So varying
 $\varpi_F$, hence $\beta$, gives other simple cuspidal representations of $G$, and our results apply equally to them. 

However varying $\varpi_F$ does not give all the simple cuspidal representations attached to the more general affine generic characters of~\S\ref{notation}. Let us analyze the situation.
We first note that an arbitrary Iwahori subgroup $I$  of  $G$   is conjugate in $G$  to our fixed Iwahori subgroup $I_{2N}$, and that its subgroups $I(1)$ and $I(2)$ are sent onto
$ I_{2N}(1)$ and $I_{2N}(2)$  by a conjugation sending $ I$  to $I_{2N}$: 
indeed $I=G_{x,0}$ is the parahoric subgroup attached to the barycenter $x$ of an alcove of the Bruhat--Tits building of $G$,  $I(1) $  is the Moy--Prasad subgroup $G_{x,0+}$  and $I(2)$ the Moy--Prasad subgroup $G_{x,(\frac{1}{2N})+}$. So we don’t get more simple cuspidal representations by choosing an Iwahori subgroup other than $I_{2N}$, and we may restrict to the ones attached to the affine generic characters of~\S\ref{notation}.

Let $\lambda$   be the affine generic character of $I_{2N}(1)$   with given parameters
$\alpha_i$   for $ i=1, \dots, N$  (which are units in $F$) and 
$\alpha_{2N}$   (which has valuation $-1$ in  $F$). Let 
$\lambda'$  be another affine generic character, with parameters $\alpha_i'$.

The same reasoning that shows that our representation $\pi$  of~\S\ref{strata} is irreducible (hence cuspidal) also shows that the intertwining of $\lambda$ and $\lambda'$  is restricted to $ZI_{2N}$, a group that normalizes $I_{2N}(1)$.
So we need to examine when $\lambda$ and $\lambda'$  are conjugate under $I_{2N}$, a result that was stated without proof in \cite[p. 21]{Oi}. We will moreover get that there are $4(q_F-1)$ isomorphism classes of simple cuspidal representations of $G $ ({\it loc. cit.}). 

Of course 
$I_{2N}(1)$   acts trivially on $\lambda$ and $\lambda'$, so it is enough to look at the conjugation action of the diagonal elements 
$d=\diag(d_1,…,d_N,1/d_N,…,1/d_1) $   in $I_{2N}$. Such an element $d$ acts on 
$\lambda$  by multiplying $\alpha_i$   (for $ i=1, \dots, N-1$)     by $d_i/d_{i+1}$, 
$\alpha_N$   by $(d_N)^2$   and $\alpha_{2N}$   by  $ (1/d_1)^2$. 
Thus conjugation by $d$ preserves the classes of $\alpha_N$      and $\alpha_{2N}$  modulo squares in $\oFx$, and also preserves 
$(\alpha_1 \cdots \alpha_{N-1})^2  \alpha_{N}\alpha_{2N}$  (which matters only modulo $1+\pF$). We easily deduce that $\lambda'$  is the conjugate of $\lambda$ by  such a diagonal element $d$ if and only if: 
\begin{enumerate}
\item   $\alpha'_{N}$   is equal to   $\alpha_{N} $   modulo squares in $\oFx$.
\item  $\alpha'_{2N}$ is equal to $\alpha_{2N}$  modulo squares in $\oFx$.
\item   $(\alpha'_1 \cdots \alpha'_{N-1})^2  \alpha'_{N}\alpha'_{2N}$ is equal to $(\alpha_1 \cdots \alpha_{N-1})^2  \alpha_{N}\alpha_{2N}$ modulo $1+\pF$.
\end{enumerate}   
(Note that given (iii), (i) is equivalent to (ii)).   

The number of conjugacy classes of $\lambda$’s is $2(q_F-1)$. Indeed let 
$\epsilon$   be a non-square in $\oFx$. Conjugating as above we may assume that 
$\alpha_i=-1$   for $ i=1, \dots, N-1$,   and that 
$\alpha_N=-1$   or $-\epsilon$, and then (iii) allows  $q_F-1$   choices for $\alpha_{2N}$. Taking the central character into account shows that indeed $G$ has $4(q_F-1)$   simple cuspidal representations up to isomorphism.

\begin{Remark} Changing the additive character $\psi$  into the character $\psi^a $  sending $x$ to $\psi(ax)$  amounts to taking $\alpha_i=a$   for $ i=1, \dots, N$   and $i=2N$. 
\end{Remark}

The results in sections  \ref{dimension1section} and \ref{dimension2Nsection} therefore apply directly to half of the simple cuspidal representations of $G$. 

To see that our results still apply to the other half, let us look at the conjugation action of $\GSp(2N,F) $    on $\Sp(2N,F)$. More precisely take the diagonal elements $d_\epsilon$ in $\GSp(2N,F) $  of the form $\diag(\epsilon,\dots,\epsilon,1,…,1) $  
where $ \epsilon$   (a non-square in $\oFx$) appears $N$  times. Then conjugation by $d_\epsilon$  preserves $I_{2N}$   and $I_{2N}(1)$, and transforms $\psi_\beta$
 into the affine generic character with 
$\alpha_i=-2$   for $ i=1, \dots, N-1$,   $\alpha_N= -\epsilon$   and $\alpha_{2N}= \frac{1}{\epsilon \varpi_F}$.   Varying $\varpi_F$   we see that we get all missing cuspidals that way.

But the reducibility points are the same for our cuspidal representation $\pi$  and its conjugate by $d_\epsilon$: indeed on $\Sp(2M+2N,F)$   we can consider the action of the similar
matrix $\diag(\epsilon,\dots,\epsilon,1,…,1)$, but this time with $N+M$ occurrences of 
$\epsilon$. Conjugating 
by that matrix on the Levi subgroup $\GL(M,F)\times \Sp(2N,F)$  induces the previous conjugation on $\Sp(2N,F)$, but \shaun{the} identity on $\GL(M,F)$.    

A consequence  of the preceding analysis is the following result, which follows from 
Propositions \ref{Jordanpi} and \ref{JordanpiGL} by conjugation inside 
$\GSp(2N,F)$: 

\oldblue{
\begin{Theorem}\label{JordanSummary} Let $\pi$ be a simple cuspidal representation of $G$, written as 
 $\pi=\cInd_{Z I_{2N}(1)}^G \chi \otimes \psi_\beta$, where 
 $\chi$ is a character of  the center $Z \simeq\{\pm 1\}$ of $G$ and $\beta^{-1}$ 
is  a uniformizer of a totally ramified extension $E$ of $F$ of degree $2N$ normalizing  $I_{2N}(1)$.\\
The Jordan set of $\pi$ is 
$
\   \Jord(\pi)=  \{(\epsilon_1, 1),  (\sigma,1)   \}   \   
$   
where  
\begin{itemize} 
\item   $\epsilon_1$ is the ramified quadratic  character of $\Fx$ characterized by 
$$\epsilon_1(   N_{E/F}(\beta))=   (-1)^{(N+1)\frac{q-1}{2}} ; $$
\item    $\sigma$ is  the simple 
 cuspidal representation of $\GL(2N,F)$  defined by 
 $$
\sigma= \cInd_{E^\times \tilde I_{2N}(1)}^{\GL(2N,F)}  \tau\otimes \psi_{2\beta}  
$$
where $\tau_{|\oEx} $ is the quadratic character of $\oEx$ and 
$$\tau(\beta) =\chi(-1)\delta(2)  \xi(\delta, \psi ).$$
\end{itemize}
\end{Theorem}
}

\subsection{A remark on epsilon factors}\label{epsilon}

For use in the next paragraph, let us remark about the $\varepsilon$-factor at
$ s=1/2$  of $\epsilon_1$  and of $\sigma$   in Theorem \ref{JordanSummary} above.
Since $\epsilon_1$ is quadratic (equal to $\delta$) on restriction to $\oFx$, we have
$\varepsilon(\epsilon_1, \frac12, \psi)=  \xi(\delta, \psi)$.  
On the other hand the factor $\varepsilon(\sigma, \frac12, \psi)$     is computed in 
\cite[Lemma 2.2]{BH}  and is equal to 
$\frac{1}{\tau(2\beta)}$  (remarking that the trace of the matrix $\beta$  is $0$). But by Proposition  \ref{JordanpiGL} we have $\tau(2\beta)= \chi(-1)  \xi(\delta, \psi)$,  so
$\varepsilon(\epsilon_1, \frac12, \psi) \varepsilon(\sigma, \frac12, \psi)  =  \chi(-1) $.  

\section{Langlands parameters for simple cuspidals}\label{Langlands}

\subsection{The characteristic zero case}\label{LanglandsA}
Let us now assume that $F$ has characteristic $0$. In that case the local Langlands correspondence has been established by Arthur, and our results about reducibility points allow us to give the parameter of a simple cuspidal representation of $G$, thus completing, in the special case of simple cuspidal representations, the results of \cite{BHS}.

\begin{Theorem}\label{LanglandsTh}  Let $\pi$  be a simple cuspidal representation of $\Sp(2N,F)$ as in Theorem~\ref{JordanSummary}. Then the parameter of $\pi$ is the direct sum of the quadratic character $\omega$ of $W_F$ corresponding to $\epsilon_1$ and an irreducible orthogonal representation of dimension $2N$, corresponding via the local Langlands correspondence for $\GL(2N,F)$  to the cuspidal representation $\sigma$ of Proposition~\ref{JordanpiGL}.  
\end{Theorem}

\begin{Remark} 
Once a local Langlands correspondence for $G$ is established when $F $  has characteristic $p$, we get the result in that case too.
There has been recent progress on establishing this correspondence when $F$ has characteristic $p$   
(see Ganapathy--Varma \cite{GV}, Gan--Lomel\'\i$\ $ \cite{GL}, and current work of Aubert and Varma). Besides, for a generic cuspidal
representation $\pi$  of $G$ (in particular for a simple cuspidal one),   Lomel\'\i$\ $ \cite{Lom}   
has used converse theorems to produce
a parameter for $\pi$. At another occasion we shall show that the arguments of the present section still apply to
explicit the parameter, giving the exact same statement.  
\end{Remark}

\begin{Remark}  Conjugating $\pi$  inside $\GSp(2N,F)$  by $d_\epsilon$  as in \ref{othersimple} gives a representation with the same parameter.
\end{Remark}

%{\color{red}
%(Ce n’est pas vraiment un bel énoncé, faut-il répéter les formules ? Rédiger différemment?) }  

\subsection{An alternative proof: method}\label{LanglandsB}
In fact, the analysis and results of \cite{BHS}, supplemented by an identity due to Lapid, are enough to get the previous theorem, without using the computations of sections \ref{dimension1section} and \ref{dimension2Nsection}, as we show presently. That gives a consistency check on those very computations, when $F$ has  characteristic $0$.

Let $\pi$  be our simple cuspidal representation as in~\S\ref{othersimple}. From \cite{BHS} we know already that the parameter $\rho$ of $\pi$  is the direct sum of a quadratic character $\omega$  of $W_F$  and an irreducible orthogonal representation $\tau$  of dimension $2N$, corresponding to a simple cuspidal representation $\sigma$  of $\GL(2N,F)$ constructed from the stratum attached to $2\beta$. In particular, $\tau$  has Swan exponent $1$, hence is not tame, and has trivial stabilizer under
character twists. In principle the results of \cite{BHS} allow us to determine the restrictions of $\omega$ and $\tau$  to the inertia group, so the only ambiguity left is small: we could twist  $\omega$ and $\tau$  by unramified quadratic characters (see \cite{BHS}, section 6, in particular 6.6 Proposition) and have an equally plausible parameter after the results of \cite{BHS}.

To remove that ambiguity we note two things. The first is that, for an unramified character $\eta$ of $W_F$ of order $2$, $\tau$ and $\eta\tau$ have equal determinant, since $\dim(\tau)$ is even. So $\omega$ is determined by 
$\omega=\det(\tau)=\det( \eta\tau)$: there is no ambiguity in $\omega$. The second is that the $\varepsilon$-factor of $\tau$ is sensitive to that character twist, because $\tau$ has Swan exponent $1$ hence Artin exponent $2N+1$: we have 
$\varepsilon( \eta\tau, \frac12, \psi)=-\varepsilon(  \tau, \frac12,\psi)$.
 But the main result of Lapid \cite{La}  gives us precisely the necessary information. 
Indeed, the representation $\pi$ is generic, so its Langlands--Shahidi factors 
$\varepsilon( \pi, s,\psi)$  
  are defined. But the local Langlands correspondence for $\Sp(2N)$ preserves the $\varepsilon$-factors, in the sense that 
$\varepsilon( \pi, s,\psi)=\varepsilon( \omega, s,\psi)\varepsilon( \tau, s,\psi)$     (for that preservation, see Appendices A and B in \cite{AHKO}). Similarly by the local Langlands correspondence for $\GL(2N,F)$, we have $\varepsilon( \tau, s,\psi)=\varepsilon( \sigma, s,\psi)$.   
  The result of Lapid says that  $\varepsilon( \pi, \frac12,\psi)$    is the value 
$\chi(-1)$ of the central character of $\pi$ at $-1$. Thus we deduce 
$\varepsilon( \sigma, \frac12,\psi)\varepsilon( \omega, \frac12,\psi) =\chi(-1) $,    which resolves the ambiguity in $\rho$.

\subsection{An alternative proof: results}\label{LanglandsC} Let us identify $\omega$ and the character corresponding to it via class field theory, also written $\omega$. Let us show that $\omega$ is the character 
$\epsilon_1$  of Theorem~\ref{JordanSummary}. We first
show that $\omega$ is ramified. Indeed $\tau$ has Artin exponent $2N+1$  and the orthogonal representation $\rho$  has trivial determinant. Then $\rho$  has even Artin exponent by an old result of Serre \cite{Se}, and that implies that $\omega$ has odd Artin exponent, hence is quadratic ramified. 

The cuspidal representation $\sigma$  of $\GL(2N,F)$ has central character $\omega$ and is constructed from the affine generic character $\psi_{2\beta}$  of the subgroup $J^1=I_N(1)$   of $\GL(2N,F)$. It is induced from an extension $\theta$ of  $\psi_{2\beta}$  to its normalizer $J$  in $\GL(2N,F)$, which is the group
$ (2\beta)^\mathbb Z F^\times J^1$, and that extension is \shaun{$\omega$} on $F^\times$, so is determined by its value $a$ on $2\beta$, subject to $a^{2N}=\omega((2\beta)^{2N})$.

However $\tau$ is self-dual, which imposes a condition on $a$. The contragredient of 
$\tau$  is induced from the character $\theta^{-1}$ of $J$. Saying that $\tau$ is self-dual therefore means that $\theta^{-1}$ intertwines with $\theta$ in $G$. But the restriction of $\theta^{-1}$ to $I_N(1)$ is the affine generic character  $\psi_{-2\beta}$, so it is sent to  $\psi_{2\beta}$  by conjugation by the diagonal matrix $\diag(1,-1,1,-1,…,1,-1)$, which conjugates  $\beta$ to $-\beta$. The condition on $a$ is therefore that $\theta (-2\beta)= \frac{1}{\theta (\shaun{2\beta})}$, that is 
$a^2=\omega(-1)$. Thus a fortiori 
$\omega(\beta^{2N})=a^{2N}=\omega(-1)^N$. But, as seen in~\S\ref{resultGL1}, 
$N_{E/F}(\beta)= -\beta^{2N}$, so  $\omega(N_{E/F}(\beta))=\omega(-1)^{N-1}$.  We happily find exactly the same recipe as in Proposition~\ref{Jordanpi}, so that indeed  $\omega=\epsilon_1$.  It now also follows from~\S\ref{epsilon} that 
$\sigma$  is given by the recipe of Theorem~\ref{JordanSummary}.

\subsection{The case of non-simple cuspidals for~$\Sp(4,F)$}\label{LanglandsD} Let us briefly comment on what Lapid’s result brings to the analysis of the examples in \cite[\S6.9]{BHS}. When $N=1$, it gives supplementary information which determines the parameter of a cuspidal representation of 
$\SL(2,F)$  (of course, that case can also be deduced from the local Langlands correspondence for $\GL(2,F)$).

Let us look at the more interesting case where $N=2$. We do not consider parameters with an occurrence of $\St_3$: the corresponding packets contain non-cuspidal discrete series, they have been determined explicitly by Suzuki and Xu \cite{SX}, thus confirming guesses of the second author decades ago (Lettre aux esp\'equatrophiles).

An ambiguous case in \cite{BHS} was that of a parameter involving $3$ quadratic characters and an irreducible orthogonal representation $\rho$  of dimension $2$, induced from a quadratic ramified extension. In that case the Artin exponent of $\rho$ is odd, so choosing between $\rho$  and the other possibility $\rho'$ (the twist of 
$\rho$  by the unramified order $2$ character) is done using Lapid’s result. However when the parameter contains two ambiguous components of dimension $2$, adding Lapid’s result does not resolve all ambiguities.

%%%%%%%%%%%%%%%%%%%%%%%%%%%%%%%%%%%%%%%%%%%%%%%%%%%%%%%%
%%% Bibliography
%%%%%%%%%%%%%%%%%%%%%%%%%%%%%%%%%%%%%%%%%%%%%%%%%%%%%%%%inertial 

\def\cprime{$'$}


\begin{thebibliography}{10} 

\bibitem{AHKO}
Moshe Adrian, Guy Henniart, Eyal Kaplan, Masao Oi. \newblock  
Simple supercuspidal $L$-packets of split special orthogonal groups over dyadic fields.
\newblock https://arxiv.org/abs/2305.09076


\bibitem{Arthur}
James Arthur. \newblock
{\em The endoscopic classification of representations}, volume 61 of {\em American Mathematical Society Colloquium Publications}. 
\newblock American Mathematical Society, Providence, RI, 2013.

\bibitem{BB} Laure Blasco and  Corinne Blondel. \newblock Alg\`ebres de Hecke et s\'eries principales g\'en\'eralis\'ees de
$\Sp_4(F)$.  \newblock {\em Proc. London Math. Soc.} (3) 85 (2002), 659--685.

\bibitem{Bl2} Corinne Blondel. \newblock  Propagation de paires couvrantes dans les groupes symplectiques.  \newblock {\em 
Representation Theory} 10 (2006),   399-434. 

\bibitem{BHS}  Corinne Blondel, Guy Henniart, Shaun Stevens.     \newblock  Jordan blocks of cuspidal representations of symplectic groups. \newblock  {\em Algebra Number Theory}    12 (2018), no. 10, 2327--2386. 
 
\bibitem{BS}  Corinne Blondel and Shaun Stevens. \newblock  
Genericity of supercuspidal representations of $p$-adic $\Sp_4$.  
\newblock {\em Compos. Math.} 145 (2009), no.1, 213--246.

\bibitem{BT} 
Corinne Blondel and Geo Kam-Fai Tam.  \newblock 
  Base change for ramified unitary groups: the strongly ramified case.  \newblock  {\em	Journal f\"ur die reine und angewandte Mathematik (Crelles Journal)} vol. 2021, no. 774 (2021)127-161.    

\bibitem{BH} 
Colin~J. Bushnell and Guy  Henniart. \newblock
Langlands parameters for epipelagic representations of  {${\rm GL}({\rm n})$}.  
\newblock {\em Math. Ann.} 358 no.1-2 (2014), 433--463.

\bibitem{BK}
Colin~J. Bushnell and Philip~C. Kutzko.
\newblock {\em The admissible dual of {${\rm GL}({\rm N})$} via compact open
  subgroups}, volume 129 of {\em Annals of Mathematics Studies}.
\newblock Princeton University Press, Princeton, NJ, 1993. 

\bibitem{BK1}
Colin~J. Bushnell and Philip~C. Kutzko.
\newblock Smooth representations of reductive $p$-adic groups: structure theory
  via types.
\newblock {\em Proc. London Math. Soc. (3)}, 77(3) (1998), 582--634.

\bibitem{BK2}
Colin~J. Bushnell and Philip~C. Kutzko.
\newblock Semisimple types in {${\rm GL}\sb n$}.
\newblock {\em Compositio Math.}, 119(1) (1999), 53--97.

\bibitem{Ca} W. Casselman. \newblock {\em Introduction to the theory of admissible
representations of $p$-adic reductive groups}, manuscript, c.1974. 
\newblock available at: 
{\it www.math.ubc.ca/people/faculty/cass/research/p-adic-book.dvi}

\bibitem{GL} 
Wee Teck Gan and Luis Lomel\'\i. \newblock 
Globalization of supercuspidal representations over function fields and applications. 
\newblock  {\em
J. Eur. Math. Soc.} 20 (2018), no.11, 2813--2858

\bibitem{GV}  
Radhika Ganapathy  and  Sandeep Varma. \newblock  
On the local Langlands correspondence for split classical groups over local function fields. \newblock {\em 
J. Inst. Math. Jussieu} 16 (2017), no.5, 987--1074.

\bibitem{GR}
Benedict~H. Gross and Mark Reeder.
\newblock Arithmetic invariants of discrete {L}anglands parameters.
\newblock {\em Duke Math. J.}, 154(3) (2010), 431--508.

\bibitem{La} Erez M.  Lapid. \newblock On the root number of representations of orthogonal type. \newblock {\em Compos. Math.} 140 (2004), no. 2, 274--287.

\bibitem{Lom} Luis Lomel\'\i. \newblock
Functoriality for the classical groups over function fields. 
\newblock {\em Int. Math. Res. Notices} (2009), 4271--4335.

\bibitem{LS}
 Jaime Lust and Shaun Stevens. \newblock 
On depth zero $L$-packets for classical groups. \newblock
{\em Proc. Lond. Math. Soc.} (3)121(2020), no.5, 1083--1120.

\bibitem{MS}
Michitaka Miyauchi and Shaun Stevens.
\newblock Semisimple types for {$p$}-adic classical groups.
\newblock {\em Math. Ann.}, 358(1-2) (2014), 257--288.

\bibitem{Moeglin}
Colette M{\oe}glin.
\newblock Paquets stables des s\'eries discr\`etes accessibles par endoscopie tordue; leur param\`etre de {L}anglands.
\newblock
pp. 295--336 in {\em Automorphic forms and related geometry: assessing the legacy of I. I. Piatetski-Shapiro (New Haven, CT, 2012)}.
\newblock Contemp. Math. 614, Amer. Math. Soc., Providence, RI, 2014.

\bibitem{Oi} Masao Oi.    
\newblock Simple supercuspidal $L$-packets of quasi-split classical groups. 
\newblock  	
https://arxiv.org/abs/1805.01400  {\it to appear in Mem. AMS}.

\bibitem{Sa2008}  Gordan Savin. \newblock  
Lifting of generic depth zero representations of classical groups.
\newblock  {\em 
J. Algebra} 319(2008), no.8, 3244--3258.   

\bibitem{Se} Jean-Pierre   Serre.  \newblock 
Conducteurs d'Artin des caract\`eres r\'eels.  \newblock  
{\em Invent. Math.} 14(1971), 173--183.

\bibitem{Sh} Freydoon Shahidi. \newblock  The notion of norm and the representation theory of  orthogonal groups.   \newblock 
{\em Invent. Math.} 119 (1995), 1--36. 

\bibitem{S1}
Shaun Stevens.
\newblock Intertwining and supercuspidal types for {$p$}-adic classical groups.
\newblock {\em Proc. London Math. Soc. (3)}, 83(1) (2001), 120--140. 

\bibitem{S5}
Shaun Stevens.
\newblock The supercuspidal representations of {$p$}-adic classical groups.
\newblock {\em Invent. Math.}, 172(2) (2008), 289--352. 

\bibitem{SX} Kenta Suzuki and  Yujie Xu. \newblock  
The explicit Local Langlands Correspondence for $\GSp_4$, $\Sp_4$ and stability (with an application to Modularity Lifting). \newblock 
https://arxiv.org/abs/2304.02622

\end{thebibliography}
\end{document}